\input amstex
\magnification=\magstep1
\documentstyle{amsppt}
\vsize 45pc
\NoBlackBoxes
\TagsAsMath
\define\a{\frak a}
\define\bb{\frak b}
\define\k{\kappa}
\define\m{\frak m}
\define\M{\frak M}

\define\p{\frak p}
\define\q{\frak q}

\define\F{\Bbb F}
\define\J{\Cal J}
\define\OO{\Cal O}

\define\N{\Bbb N}
\define\Q{\Bbb Q}
\define\R{\Bbb R}
\define\Z{\Bbb Z}
\define\Rees{{\bold R}}
\define\red{\text{red}}
\define\Ann{\operatorname{Ann}}
\define\cl{\operatorname{cl}}
\define\Cl{\operatorname{Cl}}
\define\Ker{\operatorname{Ker}}
\define\Spec{\operatorname{Spec }}
\define\Proj{\operatorname{Proj }}

\define\Supp{\operatorname{Supp}}
\define\Hom{\operatorname{Hom}}
\define\Div{\operatorname{div}}
\define\Int{\operatorname{Int}}
\define\dcap{\displaystyle\bigcap}
\define\rd#1{{\lfloor #1 \rfloor}}
\define\rup#1{{\lceil #1 \rceil}}
\define\enddemoo{\hfill $\square$ \enddemo}
\define\lra{\longrightarrow}
\NoRunningHeads

\topmatter
\title A generalization of tight closure \\ and multiplier ideals
\endtitle
\author Nobuo Hara and Ken-ichi Yoshida
\endauthor
\thanks Both authors are partially supported by Grant-in-Aid for Scientific
Research, Japan.
\endthanks
\abstract
We introduce a new variant of tight closure associated to any fixed ideal
$\a$, which we call $\a$-tight closure, and study various properties thereof.
In our theory, the annihilator ideal $\tau(\a)$ of all $\a$-tight closure
relations,
which is a generalization of the test ideal in the usual tight closure theory,
plays a particularly important role. We prove the correspondence of the
ideal $\tau(\a)$ and the multiplier ideal associated to $\a$ (or, the
adjoint of $\a$ in
Lipman's sense) in normal $\Q$-Gorenstein rings reduced from characteristic
zero to characteristic $p \gg 0$. Also, in fixed prime characteristic, we
establish some properties of $\tau(\a)$ similar to those of multiplier ideals
(e.g., a Brian\c con-Skoda type theorem, subadditivity, etc.) with considerably
simple proofs, and study the relationship between the ideal $\tau(\a)$ and the
F-rationality of Rees algebras.
\endabstract
\endtopmatter

\document

\head Introduction \endhead

The notion of tight closure, introduced by Hochster and Huneke \cite{HH1}
more than a decade ago, has emerged as a powerful new tool in commutative
algebra. Tight closure gives remarkably simple characteristic $p$ proofs
of several results that were not thought to be particularly related, e.g.,
that rings of invariants of linearly reductive groups acting on regular
rings are Cohen--Macaulay, that the integral closure of the $n$th power
of an $n$ generator ideal of a regular ring is contained in the ideal
(the Brian\c con--Skoda theorem),
and so on. Also, the notions of F-regular and F-rational rings are
defined via tight closure, and they turned out to correspond to log
terminal and rational singularities, respectively (\cite{Ha1}, \cite{HW},
\cite{MS}, \cite{Sm1}). This result is generalized to the correspondence
of test ideals and multiplier ideals (\cite{Ha2}, \cite{Sm2}), both of
which play very important roles in the tight closure theory and birational
algebraic geometry, respectively.

The test ideal of a ring $R$ of characteristic $p$, denoted by
$\tau(R)$, is the annihilator ideal of all tight closure relations
of $R$. On the other hand, the notion of multiplier ideals has a
few variants. Originally, a multiplier ideal is defined analytically
for a given plurisubharmonic function on a complex analytic manifold
\cite{N}.
This is reformulated in the algebro-geometric setting (in characteristic
zero) in terms of resolution of singularities and discrepancy divisors
(\cite{Ei}, \cite{La}). Actually, two types of 
multiplier ideals are defined in this setting, that is, the multiplier
ideal $\J(D)$ associated to a $\Q$-divisor $D$ and the one $\J(\a)$
associated to an ideal $\a$. The latter is also defined by Lipman
\cite{Li} in a more algebraic context and called the "adjoint ideal."

Precisely saying, the multiplier ideal which is proved to correspond
to the test ideal $\tau(R)$ is the one associated to the trivial divisor
$D = 0$ or the unit ideal $\a = R$, which defines the non-log-terminal
locus of $\Spec R$. In most applications, however, the usefulness
of multiplier ideals is performed by considering multiplier ideals
associated to various divisors or ideals; see e.g., \cite{Ei}, \cite{La},
\cite{Li}. Thus we are tempted to define a sort of a tight closure
operation and a "test ideal" associated to given any $\Q$-divisor or
ideal.

In this paper, we introduce a generalization of tight closure, which we
call {\it $\a$-tight closure}, associated to an ideal $\a$, study various
properties (including the relationship with multiplier ideals), and give
some applications of $\a$-tight closure. (In \cite{T}, the reader can find
an attempt to generalize tight closure to the other direction, that is,
$\Delta$-tight closure associated to a $\Q$-divisor $\Delta$; see also
\cite{HW}.)  Actually, given an ideal $\a$ of a Noetherian ring $R$ of
characteristic $p > 0$, we define the $\a$-tight closure $I^{*\a}$ of an
ideal $I \subseteq R$ to be the ideal consisting of all elements $z\in R$
for which there exists an element $c\in R$ not in any minimal prime ideal
such that
$$
cz^{p^e}\a^{p^e} \subseteq I^{[p^e]}
$$
for all $e\gg 0$, where $I^{[p^e]}$ is 
$I$ the ideal generated by the $p^e$th powers of elements of $I$; see
Definition 1.1 and also Definition 6.1 for further generalization to
"rational coefficients." We then define
the ideal $\tau(\a)$ of $R$ to be the unique largest ideal such that $\tau
(\a)I^{*\a} \subseteq I$ for all ideals $I \subseteq R$. So, in the case
where $\a =R$ is the unit ideal, the $\a$-tight closure $I^{*\a}= I^{*R}$
is equal to the tight closure $I^*$ in the usual sense, and the ideal
$\tau(\a) = \tau(R)$ is nothing but the test ideal \cite{HH1}.

There are many similarities between usual tight closure and $\a$-tight
closure, e.g., the existence of $\a$-test elements (Definition 1.6) is
proved quite similarly as in the case for usual test elements. But there
does exist a difference as well: We require a "closure" operation to
satisfy the property that, once the operation is performed, the obtained
closure does not change if one performs the operation twice or more, and
tight closure satisfies this property, namely, $(I^*)^* = I^*$. However,
it happens that $(I^{*\a})^{*\a}$ is strictly larger than $I^{*\a}$, and
so, $\a$-tight closure is not a closure operation in fact. Similarly,
unlike the usual test ideal $\tau(R)$, the ideal $\tau(\a)$ is no longer
equal to the one generated by $\a$-test elements if $\a \subsetneq R$.

In spite of the apparent disadvantage mentioned above, we find very
more advantages in the circle of ideas involving $\a$-tight closure.
The significance of the ideal $\tau(\a)$ is witnessed by the following
theorem, which ensures the expected correspondence of $\tau(\a)$ and the
multiplier ideal $\J(\a)$; see also Theorem 6.7.

\proclaim{Theorem 3.4}
Let $R$ be a normal $\Q$-Gorenstein local ring essentially of finite type
over a field and let $\a$ be a nonzero ideal. Assume that $\a \subseteq R$
is reduced from characteristic zero to characteristic $p \gg 0$, together
with a log resolution of singularities $f \colon X \to Y = \Spec R$ such
that $\a\OO_X = \OO_X(-Z)$ is invertible. Then
$$
\tau(\a) = H^0(X,\OO_X(\rup{K_X-f^*K_Y}-Z)).
$$
\endproclaim

Note that, by definition, the multiplier ideal $\J(\a)$ in
characteristic zero takes just the same form as the right-hand side
of the above equality, so, one can think of the right-hand side as a
reduction modulo $p$ of the multiplier ideal in characteristic zero.
This theorem generalizes the main results of \cite{Ha2} and \cite{Sm2},
and is proved essentially in a similar way as the preceding works,
with a little more effort.

The usefulness of multiplier ideals is tied up with vanishing theorems
in characteristic zero such as the Nadel vanishing theorem \cite{N},
which fail in characteristic $p > 0$. For example, Lipman used these
tools to establish an improved version of the Brian\c con--Skoda theorem,
which asserts that, for any ideal $\a$ of a regular local ring generated
by $r$ elements, one has $\J(\a^{n+r-1})\subseteq \a^n$ for all $n\ge 0$
\cite{Li}; see also \cite{La}. We will prove in Theorem 2.1 that the
corresponding statement $\tau(\a^{n+r-1}) \subseteq \a^n$ holds true
in characteristic $p > 0$. Taking into account the correspondence in
Theorem 3.4, we see that this gives another proof of Lipman's result.
An advantage of our prime characteristic proof here is that it is quite
elementary and simple (like the original tight closure proof of Brian\c
con--Skoda \cite{HH1}), and does not depend on desingularization or
vanishing theorems.

We shall take a look at the organization of this paper, which we hope
gives further confirmation of the usefulness of $\a$-tight closure.

After discussing basic properties of $\a$-tight closure and the ideal
$\tau(\a)$ in Section~1, we give three fundamental applications of
$\a$-tight closure in Section 2. The first one is the modified Brian\c
con--Skoda theorem via $\a$-tight closure mentioned above. Second, we
study the relationship between $\a$-tight closure and tight integral
closure defined by Hochster \cite{Ho2}, and rephrase the F-rationality
criterion of Rees algebras obtained in \cite{HWY1} in terms of $\a$-tight
closure; see Theorem 2.7 and Corollary 2.9. This enables us to give
an interesting characterization of regular local rings as the third
application. Namely, we prove in Theorem 2.15 that the regularity of
a $d$-dimensional local ring $(R,\m)$ is characterized by the property
that $\tau(\m^{d-1}) = R$. This is considered an analog of the fact
that the weak F-regularity of $R$ is characterized by the property
that $\tau(R) = R$.

Section 3 is devoted to proving the equality $\tau(\a) = \J(\a)$
in Theorem 3.4, which holds true in the situation reduced from
characteristic zero to characteristic $p \gg 0$. We note that the
containment $\tau(\a) \subseteq \J(\a)$ essentially holds true in
any fixed characteristic $p > 0$; cf.\ Proposition 3.8.

In Section 4, we establish properties of ideals $\tau(\a)$ in fixed
characteristic $p > 0$ similar to those of multiplier ideals $\J(\a)$
in characteristic zero, namely, the restriction theorem (Theorem 4.1),
the subadditivity (Theorem 4.5), and a description of ideals $\tau(\a)$
in toric case (Theorem 4.8); see \cite{DEL}, \cite{How}, \cite{La} for
the results proved for multiplier ideals. Again in light of Theorem
3.4, we can also say that the results in this section give new prime
characteristic proofs to the geometric statements for multiplier ideals
in characteristic zero, although we work with the Frobenius map in
fixed characteristic $p > 0$.

In Section 5, we explore the behavior of the ideal $\tau(I)$ for
an $\m$-primary ideal $I$ of a Gorenstein local domain $(R,\m)$ of
characteristic $p > 0$, from the viewpoint of the Rees algebra $\Rees(I)
= R[It]$ via Theorem 2.7. As a main result of this section, we prove in
Theorem 5.1 that if $\Rees(I)$ is F-rational, then its graded canonical
module is described as $\omega_{\Rees(I)} = \bigoplus_{n\ge 1} \tau(I^n)
t^n$. In particular, the equality $\tau(I) = \J(I)$ holds if $\Rees(I)$
is F-rational, and the converse is also true in case of a two-dimensional
rational double point. (It should be noted that the ideals $\tau(I)$
and $\J(I)$ may disagree in fixed positive characteristic.) Comparing
Theorem 5.1 with Hyry's results (\cite{Hy1}, \cite{Hy2}), we can also
deduce various results for $\tau(I)$.

In Section 6, we extend the notions of $\a$-tight closure
and the ideal $\tau(\a)$ to those with "rational coefficients," and
generalize some results discussed in the previous sections to the case of
rational coefficients. Although we have no explicit applications of
this generalization at the moment, we include this section for the
future reference, because recent applications of multiplier ideals
involve rational coefficients successfully (\cite{ELS}, \cite{La}).


\bigskip
Before going ahead we review part of the notation and basic notions of
the tight closure theory. We keep it minimal to avoid overlap with
Section 1. The reader is referred to Hochster and Huneke \cite{HH1--4}
and Huneke \cite{Hu} for the full development of the theory.

\subhead Notation and basic notions \endsubhead
Throughout this paper all rings are Noetherian commutative rings with unity.
For a ring $R$, we denote by $R^{\circ}$ the set of elements of $R$ which
are not in any minimal prime ideal. We will often work over a field of
characteristic $p > 0$. In this case we always use the letter $q$ for a
power $p^e$ of $p$.

Let $R$ be a Noetherian ring of characteristic $p > 0$. For an ideal $I$
of $R$ and a power $q$ of $p$, we denote by $I^{[q]}$ the ideal generated
by the $q$th powers of elements of $I$. The {\it tight closure} $I^*$ of
$I$ is defined to be the ideal consisting of all elements $z \in R$ for
which there exists an element $c \in R^{\circ}$
such that $cz^{q} \in I^{[q]}$ for all large $q = p^e$. (Tight closure
is also defined for a submodule of a module; cf.\ Section 1.)

We say that $R$ is {\it weakly F-regular} if every ideal $I$ of $R$ is
tightly closed, that is, $I^* = I$. A local ring $(R,\m)$ is said to be
{\it F-rational} if every
ideal generated by a system of parameters of $R$ is tightly closed.
In general, we say that $R$ is {\it F-rational} (resp.\ {\it F-regular})
if all of its local rings are F-rational (resp.\ weakly F-regular).

Let $F \colon R \to R$ be the Frobenius map, that is, the ring homomorphism
sending $z \in R$ to $z^p \in R$. The ring $R$ viewed as an $R$-module via
the $e$-times iterated Frobenius map $F^e \colon R \to R$ is denoted by
${}^e\! R$. We say that $R$ is {\it F-finite} if ${}^1\! R$ is a finitely
generated $R$-module. If $R$ is reduced, then $F^e \colon R \to {}^e\! R$
is identified with the natural inclusion map $R \hookrightarrow R^{1/q}$.
An F-finite reduced ring $R$ is said to be {\it strongly F-regular} if for
every element $c \in R^{\circ}$, there exists a power $q = p^e$ such that
the inclusion map $c^{1/q}R \hookrightarrow R^{1/q}$ splits as an $R$-module
homomorphism.


\head 1. Definition and basic properties of $\a$-tight closure \endhead

Let $R$ be a Noetherian ring of characteristic $p > 0$ and let $M$ be
an $R$-module. For each $e \in \N$, we denote $\F^e(M) = \F_R^e(M) :=
{}^e\! R \otimes_R M$ and regard it as an $R$-module by the action of
$R = {}^e\! R$ from the left. Then we have the induced $e$-times iterated
Frobenius map $F^e \colon M \to \F^e(M)$. The image of $z \in M$ via this
map is denoted by $z^q:= F^e(z) \in \F^e(M)$. For an $R$-submodule $N$ of
$M$, We denote by $N^{[q]}_M$ the image of the induced map $\F^e(N) \to
\F^e(M)$.

\definition{Definition 1.1}
Let $\a$ be an ideal of a Noetherian ring $R$ of characteristic $p > 0$ and
let $N \subseteq M$ be $R$-modules. The {\it $\a$-tight closure} of $N$ in
$M$, denoted by $N^{*\a}_M$, is defined to be the submodule of $M$ consisting
of all elements $z \in M$ for which there exists $c \in R^{\circ}$ such that
$$cz^q\a^q \subseteq N^{[q]}_M$$
for all large $q = p^e$. The $\a$-tight closure of an ideal $I \subseteq R$
is just defined by $I^{*\a} = I^{*\a}_R$.
\enddefinition

\demo{Remark $1.2$}
(1) In the case where $\a = R$ is the unit ideal, the $\a$-tight closure
$N^{*\a}_M = N^{*R}_M$ is nothing but the (usual) tight closure $N^*_M$
as defined in \cite{HH1}. However, unlike the usual tight closure, it may
happen that $(N^{*\a}_M)^{*\a}_M$ is strictly larger than $N^{*\a}_M$; see
Remark 1.4 (1). In this sense $\a$-tight closure is not an "honest" closure
operation in general.

(2) It seems significant to extend the definition to "rational coefficients,"
if we take into account the relationship with multiplier ideals; see
\cite{La}, \cite{T} and Sections 3 and 4. Namely, given nonnegative
$t \in \Q$ and $\a \in R$, $N \subseteq M$ as in Definition 1.1, we
can define the {\it $t\cdot\a$-tight closure} $N^{*t\cdot\a}_M$ of $N$ in $M$
and generalize some results for $\a$-tight closure to those for
$t\cdot\a$-tight closure. We treat this issue in Section 6.
\enddemo

We collect some basic properties of $\a$-tight closure in the following.
The proofs are easy and left to the reader.

\proclaim{Proposition 1.3}
Let $\a$ and $\bb$ denote ideals of a Noetherian ring $R$ of characteristic
$p > 0$ and let $L$ and $N$ denote submodules of an $R$-module $M$.
\roster
\item $N \subseteq N^{*\a}_M$ and $N^{*\a}_M/N \cong 0^{*\a}_{M/N}$.
\item If $L \subseteq N$, then $L^{*\a}_M \subseteq N^{*\a}_M$.
\item $N^{*\a\bb}_M \subseteq (N^{*\a}_M:\bb)_M$. Moreover,
if $\bb$ is a principal ideal, then the equality
$N^{*\a\bb}_M = (N^{*\a}_M:\bb)_M$ holds.
\item If $\bb \subseteq \a$, then $N^{*\a}_M \subseteq N^{*\bb}_M$. Moreover,
if $\a \cap R^{\circ} \ne \emptyset$ and if $\bb$ is a reduction of $\a$,
then the equality $N^{*\a}_M = N^{*\bb}_M$ holds.
\endroster
\endproclaim

\demo{Remark $1.4$} (1)
It follows that $N^*_M \subseteq N^{*\a}_M \subseteq (N^*_M:\a)_M$ from (3)
and (4) of Proposition 1.3. If $\a$ is a principal ideal, then the equality
on the right occurs, and $(N^{*\a}_M)^{*\a}_M = (N^*_M:\a^2)_M$ is strictly
larger than $N^{*\a}_M = (N^*_M:\a)_M$ in general.

(2) The colon capturing property \cite{HH1, Section 7} says that parameters
behave like a regular sequence modulo tight closure. Namely, if $x_1,\dots ,
x_{i+1} \in R$ are parameters, then under a mild assumption, $(x_1,\dots ,x_i)
:_R x_{i+1} \subseteq (x_1,\dots ,x_i)^*$. Since $\a$-tight closure contains
the usual tight closure, this remains true if we replace the usual tight
closure by $\a$-tight closure. In Proposition 1.5 below we slightly improve
this colon capturing property for $\a$-tight closure using the existence of a
test element. See \cite{HH1}, \cite{HH2} for the definition and detailed study
of test elements, and see also Definition 1.6 for a generalization to the
notion of $\a$-test elements.
\enddemo

\proclaim{Proposition 1.5}
Let $R$ be an equidimensional reduced excellent ring of characteristic $p > 0$
and let $\a$ be an ideal. Then for any parameters $x_1,\dots ,x_n$ in $R$,
$$
(x_1, \dots ,x_{n-1})^{*\a}:x_n \subseteq (x_1, \dots ,x_{n-1})^{*\a}.
$$
\endproclaim

\demo{Proof}
Actually, we can prove even more. Namely, let $I, J \subseteq R$ be ideals
generated by monomials in parameters $x_1,\dots ,x_n$ and $K \subseteq R$ be
the "expected" answer for $I:J$, that is, the monomial ideal which were equal
to $I:J$ if $x_1,\dots ,x_n$ formed a regular sequence. (Note that $I = K =
(x_1, \dots ,x_{n-1})$ and $J = (x_n)$ in our case.) Then we will show that
$$I^{*\a}:J \subseteq K^{*\a}.$$
Let $z \in I^{*\a}:J$. Then there is a $c \in R^{\circ}$ such that $cz^q\a^q
\subseteq I^{[q]}:J^{[q]}$ for $q = p^e \gg 0$. Since the "expected" answer
for $I^{[q]}:J^{[q]}$ is $K^{[q]}$, we have that $I^{[q]}:J^{[q]} \subseteq
(K^{[q]})^*$ by colon capturing of the usual tight closure \cite{HH1}. So,
for a test element $d \in R^{\circ}$, one has that $(cd)z^q\a^q \subseteq
K^{[q]}$ for $q = p^e \gg 0$, which means that $z \in K^{*\a}$ as required.
$\phantom{qed}$ \enddemoo

\definition{Definition 1.6}
Let $R$ be a Noetherian ring of characteristic $p > 0$ and let $\a$ be an
ideal of $R$. We say that an element $c \in R^{\circ}$ is an {\it $\a$-test
element} if $cz^q\a^q \subseteq I^{[q]}$ for all $q = p^e$ whenever $z \in
I^{*\a}$.
\enddefinition

In the case where $\a = R$ is the unit ideal, an $\a$-test element is
nothing but a test element in the usual sense \cite{HH1}. In \cite{HH2}
it is proved that a test element exists in nearly every ring of interest,
for example, in excellent reduced local rings \cite{HH2, Theorem 6.1}. We
can show that an $\a$-test element also does.

\proclaim{Theorem 1.7}
Let $R$ be a reduced Noetherian ring of characteristic $p > 0$, let $c \in
R^{\circ}$ and assume that one of the following conditions holds$:$
\roster
\item $R$ is F-finite and the localized ring $R_c$ is strongly F-regular.
\item $R$ is an excellent local ring and $R_c$ is Gorenstein and F-regular.
\endroster
Then some power $c^n$ of $c$ is an $\a$-test element for all ideals $\a
\subseteq R$.
\endproclaim

Here we prove the above theorem under assumption (1) only, according to the
method of \cite{HH0}. The case of assumption (2) is reduced to the F-finite
case by the machinery of "$\Gamma$-construction" used in \cite{HH2}. We do
not include the argument involving this reduction process because it is
somewhat long but essentially the same as that for the usual tight closure
\cite{HH2, Section 6}.

To prove the theorem in the F-finite case we need the following lemma from
\cite{HH0, Remark 3.2}, in which it is implicit that the exponent $n$ of $c$
may be independent of the choice of $d$.

\proclaim{Lemma 1.8}
Let $R$ be an F-finite reduced Noetherian ring of characteristic $p > 0$.
If the localization $R_c$ of $R$ at an element $c\in R^{\circ}$ is strongly
F-regular, then there exists an integer $n \ge 0$ depending only on $R$ and
$c$ satisfying the following property$:$ For any $d \in R^{\circ}$, there
exist a power $q'$ of $p$ and an $R$-linear map $\phi \colon R^{1/q'} \to
R$ sending $d^{1/q'}$ to $c^n$.
\endproclaim

\demo{Proof of Theorem $1.7$ in case $(1)$}
We will show that $c^n$ in Lemma 1.8 is an $\a$-test element for every $\a
\subseteq R$. Given any ideal $I$, any $z \in I^{*\a}$ and any power $q$ of
$p$, it is enough to show that $c^nz^q\a^q \subseteq I^{[q]}$. Since $z \in
I^{*\a}$, there exists $d \in R^{\circ}$ such that $dz^Q\a^Q \subseteq I^{[Q]}$
for every $Q$. Then by Lemma 1.8, there exist $q'$ and $\phi \colon R^{1/q'}
\to R$ sending $d^{1/q'}$ to $c^n$. Since $dz^{qq'}(\a^q)^{[q']} \subseteq
dz^{qq'}\a^{qq'} \subseteq I^{[qq']}$, one has
$$d^{1/q'}z^q\a^qR^{1/q'} \subseteq I^{[q]}R^{1/q'},$$
and applying $\phi$ gives $c^nz^q\a^q \subseteq I^{[q]}$, as required.
\enddemoo

\proclaim{Proposition-Definition 1.9}
Let $R$ be a Noetherian ring of characteristic $p > 0$ and let $\a$ be an
ideal of $R$. Let $E = \bigoplus_{\m} E_R(R/\m)$, the direct sum, taken over
all maximal ideals $\m$ of $R$, of the injective envelopes of the residue
fields $R/\m$. Then the following ideals are equal to each other and we
denote it by $\tau(\a)$.

\noindent$\phantom{\text{ii}}${\rm i)}
$\dcap_M \Ann_R(0^{*\a}_M)$, where $M$ runs through all finitely generated
$R$-modules.

\noindent
$\phantom{\text{i}}${\rm ii)}$\dcap_{M\subseteq E} \!\!\! \Ann_R(0^{*\a}_M)$,
where $M$ runs through all finitely generated $R$-submodules of $E$.

\noindent
If $R$ is locally approximately Gorenstein $($e.g., if $R$ is excellent
and reduced {\rm \cite{Ho1})}, then the following ideal is also equal to
$\tau(\a)$.

\noindent
{\rm iii)} $\displaystyle \bigcap_{I\subseteq R} (I:I^{*\a})$, where $I$ runs
through all ideals of $R$.
\endproclaim

\demo{Proof}
The proof is the same as that for the usual tight closure. See \cite{HH1,
Proposition 8.23} for the equality of i) and ii), and \cite{HH1, Proposition
8.25} for the equality of ii) and iii).
\enddemoo

\remark{Remark $1.10$}
In the case where $\a = R$ is the unit ideal,
$\tau(\a) = \tau(R)$ is called the test ideal. In this case, $\tau(R) \cap
R^{\circ}$ is equal to the set of test elements of $R$, and this justifies
the name "test ideal." But the name "$\a$-test ideal" for $\tau(\a)$ is
somewhat misleading if $\a \ne R$, because $\tau(\a) \cap R^{\circ}$ is not
equal to the set of $\a$-test elements in general.
\endremark

The following basic properties of the ideal $\tau(\a)$ follow from
Proposition 1.3. See Theorem 2.1 for a generalization of the latter
half of (1).

\proclaim{Proposition 1.11}
Let $R$ be a Noetherian ring of characteristic $p > 0$ and let $\a$ and $\bb$
denote ideals of $R$. 
\roster
\item $\tau(\a)\bb \subseteq \tau(\a\bb)$.
Moreover, if $\bb$ is a principal ideal of a complete local ring,
then $\tau(\a)\bb = \tau(\a\bb)$.
\item If $\bb \subseteq \a$, then $\tau(\bb) \subseteq \tau(\a)$. Moreover,
if $\a \cap R^{\circ} \ne \emptyset$ and if $\bb$ is a reduction of $\a$, then
the equality $\tau(\bb) = \tau(\a)$ holds.
\item If $R$ admits a test element and if $\a \cap R^{\circ} \ne \emptyset$,
then $\tau(\a) \cap R^{\circ} \ne \emptyset$.
\item If $R$ is weakly F-regular, then $\a \subseteq \tau(\a)$.
Moreover, if $\a$ is an ideal of pure height one, then $\a = \tau(\a)$.
\endroster
\endproclaim

\demo{Proof}
The former half of (1) is immediate from Proposition 1.3 (3). To prove
the latter half, let $(R,\m)$ be a complete local ring and let $\bb$ be
principal. Then by the Matlis duality, $\Ann_E(\tau(\a))$ is equal to
the union of $0^{*\a}_M$ taken over all finitely generated submodules
$M$ of $E = E_R(R/\m)$. So, if $z \in \Ann_E(\tau(\a)\bb)$, there exists
a finitely generated submodule $M \subset E$ such that $z \in (0^{*\a}_M:
\bb)_E$. Replacing $M$ by $M + Rz \subset E$, one has $z \in (0^{*\a}_M:
\bb)_M = 0^{*\a\bb}_M$ by Proposition 1.3 (3). Hence
$$\tau(\a\bb) = \bigcap_{M\subseteq E} \Ann_R(0^{*\a}_M:\bb)_M
               = \Ann_R(\Ann_E(\tau(\a)\bb)) = \tau(\a)\bb.$$

(2) follows from Proposition 1.3 (4), and (3) and the former half of (4)
from $\tau(R)\a \subseteq \tau(\a)$. As for the latter half of (4), it
suffices to show the following claim, since weakly F-regular rings are
normal.
\proclaim{Claim 1.11.1}
If $R$ is normal and $\a$ is an ideal of pure height one, then
$\tau(\a) \subseteq \a$.
\endproclaim

To prove the claim, considering a primary decomposition of $\a$, we may
assume without loss of generality that $\a$ is a primary ideal such that
$\p = \sqrt{\a}$ is a height one prime ideal. Then, since $R_{\p}$ is a
discrete valuation ring, we can choose $b \in \a$ such that $\a R_{\p} =
bR_{\p}$. Then $bR : \a \subseteq \a^{*\a}$. Indeed, if $z \in bR : \a$,
then $z^q\a^q \subseteq b^qR \subseteq \a^{[q]}$ for all $q = p^e$, so
that $z \in \a^{*\a}$. It now follows from $bR:\a \not\subseteq \p$ that
$\a^{*\a} \not\subseteq \p$. Since $\a$ is $\p$-primary, we have
$\tau(\a) \subseteq \a:\a^{*\a} = \a$, as claimed.
\enddemoo
\proclaim{Proposition 1.12
{\rm (cf.\ \cite{B}, \cite{HH1, Proposition 4.12})}}
Let $R \subseteq S$ be a pure ring extension of Noetherian rings of
characteristic $p > 0$ such that $R^{\circ} \subseteq S^{\circ}$.
Then for any ideal $\a$ of $R$, one has
$\tau(\a S) \cap R \subseteq \tau(\a)$.
\endproclaim

\demo{Proof}
For a finitely generated $R$-module $M$, the natural map $M=M\otimes_R R
\to M \otimes_R S$ is injective by the purity of $R \subseteq S$. Since
$R^{\circ} \subseteq S^{\circ}$, we see easily that $0^{*\a}_M \subseteq
0^{*\a S}_{M \otimes S}$ via the inclusion map $M \hookrightarrow
M \otimes_R S$. Hence, if $c \in \tau(\a S) \cap R$, then $c$ kills
$0^{*\a}_M$ for all finitely generated $R$-module $M$, so that $c \in
\tau(\a)$.
\enddemoo

By definition, the ideal $\tau(\a)$ is the annihilator of $\a$-tight closure
relations for all ideals or finitely generated modules. It will be very useful
if $\tau(\a)$ is determined by $\a$-tight closure relations for a single ideal
or a single module. Let us take a look at some cases where this is true.

\proclaim{Theorem 1.13}
Let $(R,\m)$ be a $d$-dimensional excellent normal local ring of
characteristic $p > 0$, $\a$ an ideal of $R$ and let $J \subseteq R$ be a
divisorial ideal such that the  divisor class $\operatorname{cl}(J) \in
\operatorname{Cl}(R)$ has a finite order.
Then
$$0^{*\a}_{H_{\m}^d(J)} = \bigcup_{M\subset H_{\m}^d(J)} 0^{*\a}_M,$$
where $M$ runs through all finitely generated $R$-submodules of $H_{\m}^d(J)$.
In particular, if $R$ is $\Q$-Gorenstein, then
$$\tau(\a) = \Ann_R(0^{*\a}_E),$$
where $E = E_R(R/\m) \cong H_{\m}^d(\omega_R)$.
\endproclaim

\demo{Proof}
Again the proof is the same as that for the usual tight
closure,\footnote{
The proof in \cite{Ha2, Appendix} has a minor gap at the bottom of
p.\ 1904, although the result \cite{Ha2, Theorem 1.8} itself and the arguments
in the cited references \cite{Mc}, \cite{Wi} are valid.}
but we sketch a proof according to \cite{Sm2, Lemma 3.4}, which is
based on the idea of \cite{AM, 3.1}.

Let $r$ be the order of $\cl(J) \in \Cl(R)$ and let $J^{(r)} = x_1R$.
We may assume without loss of generality that $x_1 \in \m$. Then there
exists $x_2 \in R$ and $0 \ne a \in J$ such that $x_2J \subseteq aR$,
and $x_1,x_2$ extends to a system of parameters $x_1,x_2,\dots ,x_d$
for $R$.

The point of the proof is that $\F^e(H_{\m}^d(J)) \cong H_{\m}^d
(J^{(p^e)})$ is computed by
$$H_{\m}^d(J^{(q)}) = \lim_{\to} R/(x_1^sJ^{(q)},x_2^s,\dots ,x_d^s),$$
where the direct limit map $R/(x_1^sJ^{(q)},x_2^s,\dots ,x_d^s) \to
R/(x_1^{s+1}J^{(q)},x_2^{s+1},\dots ,x_d^{s+1})$ is the multiplication
by $x_1x_2\cdots x_d$. Then an element $\xi \in H_{\m}^d(J)$ is
represented by $z$ mod $(x_1^sJ,x_2^s,\dots ,x_d^s) \in R/(x_1^sJ,
x_2^s,\dots ,x_d^s)$ for some $z \in R$ and $s \in \N$, and $\xi = $
$[z$ mod $(x_1^sJ,x_2^s,\dots ,x_d^s)]$ is mapped to $\xi^{p^e} =
[z^{p^e}$ mod $(x_1^{p^es}J^{(p^e)},x_2^{p^es},\dots ,x_d^{p^es})]$
by the $e$-times iterated Frobenius map $F^e \colon H_{\m}^d(J) \to
H_{\m}^d(J^{(p^e)})$.

Now say that $\xi \in 0^{*\a}_{H_{\m}^d(J)}$. Then there exists $c
\in R^{\circ}$ such that $c\xi^q\alpha = [cz^q\alpha$ mod $(x_1^{qs}
J^{(q)},x_2^{qs},\dots ,x_d^{qs})] = 0$ for all $q = p^e \gg 0$ and
$\alpha \in \a^q$. Since $\a^q$ is a finitely generated ideal for each
$q = p^e$, there
exists $t_e \in \N$ such that $cz^q(x_1\cdots x_d)^{t_e}\a^q \subseteq
(x_1^{qs+t_e}J^{(q)},x_2^{qs+t_e},\dots ,x_d^{qs+t_e}) \subseteq
(x_1^{qs+t_e+\rd{q/r}},x_2^{qs+t_e},\dots ,x_d^{qs+t_e})$. Then one
has $cz^q\a^q \subseteq (x_1^{qs+\rd{q/r}},x_2^{qs},\dots ,x_d^{qs})^*$
by colon capturing. Replacing $c$ by $cc'$ with $c'$ a test element
and multiplying by $x_1$, we see that $cx_1z^q\a^q \subseteq (x_1^{qs}
J^{(q)},x_2^{qs},\dots ,x_d^{qs})$. This gives
$$cx_1(x_1\cdots x_d)^qz^q\a^q
        \subseteq (x_1^{q(s+1)}a^q,x_2^{q(s+1)},\dots ,x_d^{q(s+1)})
        \subseteq (x_1^{s+1}J,x_2^{s+1},\dots ,x_d^{s+1})^{[q]}$$
for all $q= p^e \gg 0$, whence $(x_1\cdots x_d)z \in (x_1^{s+1}J,
x_2^{s+1},\dots ,x_d^{s+1})^{*\a}$. Hence $\xi$ is in the $\a$-tight
closure of zero in the cyclic (hence finitely generated) submodule
of $H_{\m}^d(J)$ generated by the image of $R/(x_1^{s+1}J,x_2^{s+1},
\dots ,x_d^{s+1})$.
\enddemoo

\demo{Discussion $1.14$}
In Sections 2 and 5, we will consider when the equality $\tau(\a) = R$ holds.
In the case $R$ is a Gorenstein local ring, one can check this condition
looking only at the $\a$-tight closure of a single parameter ideal, as we
will see below; cf.\ \cite{FW}.

Let $(R,\m)$ be a $d$-dimensional Cohen--Macaulay local ring, $\a$ any
ideal of $R$ and let $J$ be the ideal generated by a system of parameters
$x_1,\dots ,x_d$. Then
$H^d_{\m}(R) \cong {\displaystyle \lim_{\to}} R/(x_1^t,\dots ,x_d^t)$,
and $R/J$ and $H^d_{\m}(R)$ have the same socle in common via the natural
inclusion map $R/J \hookrightarrow H^d_{\m}(R)$. Then $0^{*\a}_{H^d_{\m}(R)}
= 0$ if and only if $0^{*\a}_{R/J} = 0$, or equivalently, if $J^{*\a} = J$.
In particular, the condition that $J$ is $\a$-tightly closed does not
depend on the choice of a parameter ideal $J$.

Now assume further that $(R,\m)$ is Gorenstein. Then $E = E_R(R/\m) \cong
H^d_{\m}(R)$ and one sees easily that
$\tau(\a)= \Ann_R(0^{*\a}_{H^d_{\m}(R)})
           = \bigcap_{t\in\N}(x_1^t,\dots ,x_d^t):(x_1^t,\dots ,x_d^t)^{*\a}$.
Therefore $\tau(\a) = R$ if and only if $J^{*\a} = J$ for some (or
equivalently, every) ideal $J$ generated by a system of parameters.
\enddemo

As we have seen so far, the $\a$-tight closure of the zero submodule in
the injective envelope $E_R(R/\m)$ or the top local cohomology $H_{\m}^d
(R)$ of a local ring $(R,\m)$ plays a particularly important role. We
close this section by the following proposition, which generalizes Smith's
characterization of the usual tight closure of the zero submodule in
$H_{\m}^d(R)$.

\proclaim{Proposition 1.15 {\rm (cf.\ \cite{Sm1})}}
Let $(R,\m)$ be a $d$-dimensional excellent normal local ring of
characteristic $p > 0$ and let $\a \subseteq R$ be an ideal such that
$\a \cap R^{\circ} \ne \emptyset$.
Then $0^{*\a}_{H^d_{\m}(R)}$ is the unique maximal proper submodule $N$
with respect to the property
$$\a^q F^e(N) \subseteq N \text{ for  all } q = p^e,$$
where $F^e \colon H_{\m}^d(R) \to H_{\m}^d(R)$ is the $e$-times iterated
Frobenius induced on $H_{\m}^d(R)$.
\endproclaim

\demo{Proof}
Let $c \in R^{\circ}$ be an element such that $R_c$ is regular. Then
$\hat{R}_c$ is also regular by the excellence of $R$. Hence some power
$c^n$ of $c$ is an $\a$-test element and an $\a\hat{R}$-test element by
Theorem1.7. It is easy to see that $c^n$ also works as a test element for
both the $\a$-tight closure and the $\a\hat{R}$-tight closure of the zero
submodule in $H^d_{\m}(R) = H^d_{\m\hat{R}}(\hat{R})$ (see the proof of
Theorem 1.13). Then it follows that
$0^{*\a}_{H^d_{\m}(R)} = 0^{*\a\hat{R}}_{H^d_{\m\hat{R}}(\hat{R})}$,
so we may assume without loss of generality that $R$ is a complete
local ring.

It is easy to see that $\a^qF^e(0^{*\a}_{H^d_{\m}(R)}) \subseteq
0^{*\a}_{H^d_{\m}(R)}$ for all $q = p^e$. Also, $0^{*\a}_{H^d_{\m}(R)}$
is a proper submodule of $H_{\m}^d(R)$, because it is annihilated by
$\tau(\a)$ by Theorem 1.13 and $\tau(\a) \cap R^{\circ} \ne \emptyset$
by Proposition 1.11 (3).
To prove the maximality of $0^{*\a}_{H^d_{\m}(R)}$, suppose that
$N \subset H^d_{\m}(R)$ is a proper submodule such that $\a^q F^e(N)
\subseteq N$ for all $q = p^e$. Then the Matlis dual of the exact
sequence $0 \to N \to H^d_{\m}(R) \to H^d_{\m}(R)/N \to 0$ is
$$0 \to [H^d_{\m}(R)/N]^{\vee} \to \omega_R \to N^{\vee} \to 0,$$
where $[H^d_{\m}(R)/N]^{\vee}$ is a nonzero submodule of $\omega_R$, so
both $[H^d_{\m}(R)/N]^{\vee}$ and $\omega_R$ are torsion-free $R$-module
of rank 1. Therefore $N^{\vee}$ is a finitely generated torsion module,
so that there exists $c \in R^{\circ}$ such that $cN^{\vee} = 0$. This
implies that $cN = cN^{\vee\vee} = 0$, so that $c\a^qF^e(N) = 0$ for all
$q = p^e$. Hence $N \subseteq0^{*\a}_{H^d_{\m}(R)}$, as required.
\enddemoo


\head 2. $\a$-tight closure and its applications \endhead

In this section we give some fundamental applications of $\a$-tight closure.

\subhead
Modified Brian{\c c}on-Skoda theorem via $\a$-tight closure
\endsubhead
One of the important applications of tight closure theory \cite{HH1} is a prime
characteristic proof of the Brian{\c c}on--Skoda theorem \cite{BS}, which was
originally proved by an analytic method. Later, Lipman \cite{Li} improved this
in terms of adjoint ideals. The following is a prime characteristic analog of
Lipman's "modified Brian{\c c}on--Skoda" \cite{Li, Theorem 1.4.1}; see also
Remark 3.2 (1).

\proclaim{Theorem 2.1}
Let $R$ be a Noetherian ring of characteristic $p > 0$. If $\a \subseteq R$
is an ideal generated by $r$ elements, then
$$\tau(\a^{n+r-1}) \subseteq \a^n$$
for all $n \ge 0$. If we assume further that $R$ is a complete local
ring, then
$$\tau(\a^{n+r-1}\bb) \subseteq \tau(\bb)\a^n$$
for all $n \ge 0$ and all ideals $\bb \subseteq R$.
\endproclaim

\demo{Proof}
Let $\bb \subseteq R$ be any ideal, $M$ any finitely generated $R$-module,
and suppose that $z \in (0^{*\bb}_M:\a^n)_M$. Then there exists $c \in
R^{\circ}$ such that $cz^q(\a^n)^{[q]}\bb^q = 0$ in $\F^e(M)$ for all $q =
p^e \gg 0$. Since $\a^{(n+r-1)q} \subseteq (\a^n)^{[q]}$ by the assumption,
this implies that $cz^q(\a^{n+r-1}\bb)^q = 0$ in $\F^e(M)$ for all $q =
p^e \gg 0$, so that $z \in 0^{*\a^{n+r-1}\bb}_M$. Thus
$(0^{*\bb}_M:\a^n)_M \subseteq 0^{*\a^{n+r-1}\bb}_M$, and in particular,
$\Ann_M(\a^n) \subseteq 0^{*\a^{n+r-1}}_M$. Taking the intersection of the
annihilator ideals over all finitely generated $R$-submodules $M \subseteq
E = \bigoplus_{\m} E_R(R/\m)$, we obtain
$$\tau(\a^{n+r-1}) \subseteq \bigcap_{M\subseteq E} \Ann_R(\Ann_M(\a^n))
                              = \Ann_R(\Ann_E(\a^n)) = \a^n.$$

Now assume that $(R,\m)$ is a complete local ring.
Then $\Ann_E(\tau(\bb)) = \bigcup_{M\subseteq E} 0^{*\bb}_M$ by the Matlis
duality, and it follows as in the latter half of Proposition 1.11 (1) that
$\Ann_E(\tau(\bb)\a^n)
        \subseteq \bigcup_{M\subseteq E} (0^{*\bb}_M:\a^n)_M
        \subseteq \bigcup_{M\subseteq E} 0^{*\a^{n+r-1}\bb}_M$,
where the unions are taken over all finite generated submodules $M$ of $E$.
Thus we conclude
$$\tau(\a^{n+r-1}\bb) \subseteq
\bigcap_{M\subseteq E} \Ann_R(0^{*\bb}_M:\a^n)_M
               = \Ann_R(\Ann_E(\tau(\bb)\a^n)) = \tau(\bb)\a^n.$$
\enddemoo

\demo{Remark $2.2$}
The tight closure version of the Brian\c con--Skoda theorem \cite{HH1,
Theorem 5.4} says that if $\a$ is generated by $r$ elements, then
$\overline{\a^{n+r-1}} \subseteq (\a^n)^*$ for all $n \ge 0$, where
$\overline\bb$ denotes the integral closure of an ideal $\bb$. This
implies that $\tau(R)\overline{\a^{n+r-1}} \subseteq \a^n$, and in
the case where  the test ideal $\tau(R)$ is a strong test ideal (this
is the case if $R$ is a reduced complete local ring \cite{Vr}),
$\tau(R)\overline{\a^{n+r-1}} \subseteq \tau(R)\a^n$.
Theorem 2.1 may be considered a slight improvement of these assertions,
because
$\tau(R)\overline{\a^{n+r-1}} \subseteq \tau(\overline{\a^{n+r-1}})
                                       = \tau(\a^{n+r-1})$
by basic properties of $\a$-tight closure; see also Discussion 5.2.

Recently, using arguments similar to the as above, the first author
and S.~Takagi proved a sharpened version of Theorem 2.1 \cite{HT};
cf.\ \cite{Li}, \cite{La}: if $(R,\m)$ is a complete local ring of
characteristic $p > 0$ and if $\a$ is an ideal with a reduction
generated by $r$ elements, then
$\tau(\a^{n+r-1}) = \tau(\a^{r-1})\a^n$ for all $n \ge 0$.
\enddemo

\proclaim{Corollary 2.3}
Let $R$ be an reduced excellent ring of characteristic $p > 0$ and let
$\a$ be an ideal such that $\a \cap R^{\circ} \ne \emptyset$. Then for
any $R$-modules $N \subset M$ and any $z \in M$, the following conditions
are equivalent.
\roster
\item $z \in N^{*\a}_M$, i.e., there exists $c \in R^{\circ}$ such that
$cz^q\a^q \subseteq N^{[q]}_M$ for all $q = p^e$.
\item There exists $c \in R^{\circ}$ such that
$cz^q\tau(\a^q) \subseteq N^{[q]}_M$ for all $q = p^e$.
\item There exists $c \in R^{\circ}$ such that
$cz^q\overline{\a^q} \subseteq N^{[q]}_M$ for all $q = p^e$.
\endroster
\endproclaim

\demo{Proof}
To prove (1) $\Rightarrow$ (2), choose $d \in \a^{r-1}\cap R^{\circ}$
and apply $d\tau(\a^q) \subseteq \tau(\a^{q+r-1}) \subseteq \a^q$. As
for (2) $\Rightarrow$ (3), choose a test element $d \in \tau(R) \cap
R^{\circ}$ and note that $d\overline{\a^q} \subseteq \tau(\a^q)$.
$\phantom{qed}$
\enddemoo

\proclaim{Corollary 2.4}
Let $(R,\m)$ be a $d$-dimensional Noetherian local ring of characteristic
$p > 0$ with infinite residue field. Then for any ideal $\a \subseteq R$
and for any $n > 0$, one has
$\tau(\a^{n+d-1}) \subseteq \a^n$.
If, in addition, $(R,\m)$ is complete, then
$\tau(\a^{n+d-1}) \subseteq \tau(R)\a^n$.
\endproclaim

\demo{Proof}
First assume that $\a$ is an $\m$-primary ideal. Then $\a$ has a minimal
reduction $\q \subseteq \a$ generated by $d$ elements, so that
$$\tau(\a^{n+d-1}) = \tau(\q^{n+d-1}) \subseteq \q^n \subseteq \a^n.$$
Next let $\a$ be an arbitrary ideal. Since our assertion holds true for
every $\m$-primary ideal, it follows that
$$\tau(\a^{n+d-1}) \subseteq \bigcap_{t\in \N} \tau((\a+\m^t)^{n+d-1})
                   \subseteq \bigcap_{t\in \N} (\a+\m^t)^n = \a^n.$$
The latter half is proved in a similar way. \enddemoo

\subhead
Tight integral closure vs.~$\a$-tight closure
\endsubhead
First, we recall the notion of tight integral closure, which was
introduced by Hochster.

\definition {Definition 2.5 {\rm (\cite{Ho2})}}
Let $R$ be a Noetherian ring, and let $\{I_1,\,\ldots,I_n\}$ be a set of
ideals in $R$. An element $x \in R$ is in the {\it tight integral closure}
$\{I_1,\,\ldots,I_n\}^{\underline{*}}$ if there exists $c \in R^{\circ}$
such that $cx^q \in \sum_{i=1}^n I_i^q$ for all sufficiently large $q =
p^e$.
\enddefinition

In \cite{HWY1}, the present authors have studied the F-rationality of
Rees algebras $\Rees(I) = R[It]$ for $\m$-primary ideals $I$, jointly
with K.-i.~Watanabe. One of the main results in \cite{HWY1} is the
following theorem which gives a criterion for F-rationality of Rees
algebras in terms of tight integral closure.

\proclaim {Theorem 2.6 {\rm (cf.\ \cite{HWY1, Theorem 2.2})}}
Let $(R,\m)$ be an excellent Cohen--Macaulay normal local ring of
characteristic $p > 0$ with infinite residue field. Let $I$ be an
$\m$-primary ideal of $R$ and $J$ its minimal reduction. Fix any
system of parameters $x_1,\,\ldots,x_d$ for $R$ generating $J$ and
put $J^{[l]} = (x_1^l,\,\ldots,x_d^l)$ for $l \ge 1$. Then the Rees
algebra $\Rees(I) = R[It]$ is F-rational if and only if\/ $\Rees(I)$
is Cohen--Macaulay and the following equalities hold$:$
$$
\{I^{dl-r},\,x_1^lR,\dots,\,x_d^lR\}^{\underline *} = I^{dl-r}+J^{[l]}
              \text{ for all $l,\,r \ge 1$ with $1 \le r \le dl-1$.}$$
\endproclaim

We will show that all tight integral closures appearing in the above
theorem can be represented as the form of some \lq\lq $\a$-tight
closure." Namely, we have the following theorem.

\proclaim {Theorem 2.7}
Let $(R,\m)$ be an excellent equidimensional reduced local ring of
characteristic $p>0$ with $d = \dim R \ge 1$. Also, let $x_1,\,\ldots,
x_d$ be a system of parameters of $R$, and put $J=(x_1,\,\ldots,x_d)R$.
Then we have
$$\{J^{dl-r},\,x_1^lR,\,\ldots,x_d^lR\}^{\underline{*}}
                          =(x_1^l,\,\ldots,x_d^l)^{*J^r}$$
for all integers $l,\,r \ge 1$.
\endproclaim

Before proving the above theorem, we give some corollaries.
We now recall that
$$\overline{I_1} + \cdots + \overline{I_n} \subseteq
                    \{I_1,\,\ldots,I_n\}^{\underline{*}};$$
see \cite{Ho, Proposition 1.4}.

\proclaim {Corollary 2.8}
Under the same notation as in Theorem $2.7$,
if $J \subseteq I \subseteq \overline{J}$,
then $(J^{[l]})^{*J^r} \supseteq \overline{J^{dl-r}} +J^{[l]}$
for all $l,\,r \ge 1$. In particular, we have
\roster
\item $J + \overline{I^{d-1}} \subseteq J^{*I}$.
\item If $J^{*I^r} = J$, then $I^{d-r} \subseteq J$.
\endroster
\endproclaim

As an application of Theorem 2.7, we can rewrite Theorem 2.6 as follows.

\proclaim {Corollary 2.9}
Under the same notation as in Theorem $2.6$, the Rees
algebra $\Rees(I)$ is F-rational if and only if\/ $\Rees(I)$
is Cohen--Macaulay and the following equalities hold$:$
$$(J^{[l]})^{*I^r} = I^{dl-r}+J^{[l]}
        \text{ for all $\;l,\,r \ge 1$ with $1 \le r \le dl-1$.}$$
\endproclaim

\proclaim {Corollary 2.10}
Let $R$ be an excellent F-rational local ring with $\dim R = 2$.
Then for any parameter ideal $J$ of $R$, we have
$J^{*J} = \overline{J}$.
\endproclaim

\demo{Proof}
Put $I = \overline{J}$. Then $J$ is a minimal reduction of $I$.
Since $\Rees(I)$ is F-rational by \cite{HWY1, Theorem 3.1}, it follows
from Theorems 2.6 and 2.7 that $J^{*J} = J^{*I} = I+J =I$, as required. 
\enddemoo

In the following, we prove Theorem 2.7, and so we assume that $(R,\m
,k)$ is an excellent equidimensional (not necessarily reduced) local
ring of characteristic $p > 0$ with $d = \dim R \ge 1$. Also, let
$x_1,\,\ldots,x_d$ be a system of parameters of $R$ and put $J =
(x_1,\,\ldots,x_d)R$.
Further, we set $J^{[l]}:=(x_1^l,\,\ldots,x_d^l)$ for all $l \ge 1$.

\proclaim {Lemma 2.11}
Suppose that $x_1,\,\ldots,x_d$ form a regular sequence.
Then for all integers $l,\,r\ge 1$ we have
$$
(x_1^l,\,\ldots,x_d^l) : (x_1,\,\ldots,x_d)^r =
(x_1,\,\ldots,x_d)^{dl-r-d+1} + (x_1^l,\,\ldots,x_d^l),
$$
that is, $J^{[l]} : J^r = J^{dl-r-d+1} +J^{[l]}$.
\endproclaim

\demo{Proof}
The right-hand side is contained in the left-hand side because
$J^{dl-r-d+1}J^r = J^{d(l-1)+1} \subseteq J^{[l]}$. We must show the
opposite inclusion. To do that, let $w =x_1^{a_1}\cdots x_d^{a_d}$,
where $0 \le a_i \le l-1$ for all $i$. First suppose that
$\sum_{i=1}^d a_i = dl-r-d$. If we put $b_i = l-1-a_i$ for all $i$,
then $b_i \ge 0$,\, $\sum_{i=1}^d b_i = r$ and
$w\cdot x_1^{b_1}\cdots x_d^{b_d} = x_1^{l-1}\cdots x_d^{l-1} \notin J^{[l]}$.
Next suppose that $\sum_{i=1}^d a_i \ge dl-r-d+1$.
Then for all integers $c_i \ge 0$ with $\sum_{i=1}^d c_i =r$, we have
$w \cdot x_1^{c_1}\cdots x_d^{c_d}
                   = x_1^{a_1+c_1} \cdots x_d^{a_d+c_d} \in J^{[l]}$
because $\sum_{i=1}^d (a_i+c_i) \ge d(l-1)+1$.
Since $J^{[l]}:J^r$ is generated by monomials in $x_1,\,\ldots,x_d$,
the assertion follows from the above argument. 
\enddemoo

Using the colon capturing property of tight closure
we obtain the following 

\proclaim {Corollary 2.12}
Under the above notation, for all $l,\,r \ge 1$, we have
$$
J^{[l]}:J^r \subseteq (J^{dl-r-d+1}+J^{[l]})^{*}.
$$
\endproclaim

\demo{Proof}
First suppose that $R$ is complete, reduced. If we put
$S = k[[x_1,\,\ldots,x_d]]$, then $S$ is a complete regular local domain
and $R$ is a finitely generated torsion-free $S$-module. Also, if we put
$J_0=(x_1,\,\ldots,x_d)S$ and $J_0^{[l]} = (x_1^l,\,\ldots,x_d^l)S$,
then $J=J_0R$ and $J^{[l]} = J_0^{[l]}R$. Using the colon capturing property
of tight closure and the previous lemma, we get
$$
J^{[l]} : J^r \subseteq \left((J_0^{[l]} : J_0^{r})R\right)^{*} =
\left((J_0^{dl-r-d+1}+J_0^{[l]})R\right)^{*}
= \left(J^{dl-r-d+1}+J^{[l]}\right)^{*}.
$$
\par
Next we consider the general case. Fix $l,\,r \ge 1$ and put
$K = J^{dl-r-d+1}+J^{[l]}$. Applying the above argument to
$\widehat{R}_{\red} = \widehat{R_{\red}}$, we have
$$
J^{[l]}\widehat{R}_{\red} : J^r\widehat{R}_{\red}
\subseteq (K\widehat{R}_{\red})^{*} =(K\widehat{R})^{*}\widehat{R}_{\red}
$$
and hence
$J^{[l]}\widehat{R} : J^r\widehat{R} \subseteq (K\widehat{R})^{*}$.
By \cite{BH, Proposition 10.3.18}, we get
$$
J^{[l]}:J^r = (J^{[l]}\widehat{R} : J^r\widehat{R})\cap R
= (K\widehat{R})^{*} \cap R =K^{*}\widehat{R} \cap R =K^{*},
$$
as required. 
\enddemoo

We are now ready to prove Theorem 2.7.

\demo{Proof of Theorem $2.7$}
Note that $R$ admits a test element $c' \in R^{\circ}$ because $R$ is
an excellent reduced local ring (\cite{HH2, Theorem 6.1}).
\par
Let $z \in (J^{[l]})^{*J^r}$. By definition, there exists $c'' \in
R^{\circ}$ such that $c''z^q J^{rq} \subseteq J^{[lq]}$ for all
$q=p^e,\,e \gg 0$. Corollary 2.12 implies that
$$
c''z^q \in J^{[lq]}:J^{rq}
\subseteq (J^{(dl-r)q-d+1}+J^{[lq]})^*
$$
and hence
$$
c'c''z^q \in J^{(dl-r)q-d+1} + J^{[lq]}
$$
for all $q=p^e,\,e \gg 0$. Take any element $c''' \in J^{d-1} \cap
R^{\circ}$ and put $c = c'c''c''' \in R^{\circ}$.
Then $cz^q \in J^{(dl-r)q}+J^{[lq]}$ for all $q=p^e,\,e \gg 0$.
Thus $z \in \{J^{dl-r},\,x_1^lR,\,\ldots,x_d^lR\}^{\underline{*}}$.

Next we prove the opposite inclusion.
Let $w \in \{J^{dl-r},\,x_1^lR,\,\ldots,x_d^lR\}^{\underline{*}}$.
By definition, there exists $c'' \in R^{\circ}$ such that
$c''w^q \in J^{(dl-r)q} + J^{[lq]}$ for all $q=p^e, \,e \gg 0$.
Thus $c''w^q J^{rq} \subseteq J^{dlq} +J^{[lq]}$.

On the other hand, by virtue of the tight closure Brian\c con--Skoda
\cite{HH1, Theorem 5.4},
we have
$$
J^{dlq} \subseteq \overline{J^{dlq}} \subseteq (J^{[lq]})^{*}.
$$
Taking a test element $c' \in R^{\circ}$, we have
$c'\left(J^{[lq]}\right)^{*} \subseteq J^{[lq]}$ for all $q=p^e,\,e \gg 0$.
In particular, we have $c'c''w^q J^{rq} \subseteq J^{[lq]}$ for all
sufficiently large $q=p^e$ and thus $w \in (J^{[l]})^{*J^r}$, as required.
\enddemoo

\remark{Remark $2.13$}
Although we can easily see that
$J^{*I}R_{\red} \subseteq (JR_{\red})^{*IR_{\red}}$ always holds,
we do not know whether or not the opposite inclusion holds.
If it were true, we could remove the assumption that \lq\lq $R$ is reduced''
in Theorem $2.7$.
\endremark

\subhead
A characterization of regular local rings \endsubhead
Let $(R,\m,k)$ be an excellent equidimensional reduced local ring of
characteristic $p>0$ with $d =\dim R \ge 1$.
Then $R$ is F-rational if and only if $J^{*} = J$ for some (every)
parameter ideal $J$ of $R$; see \cite{HH2} and \cite{FW}.

Let $J$ be a minimal reduction of $\m$. Since $J^{*} = J^{*R}$, we have
the following increasing sequence of ideals in $R$$:$
$$
J \subseteq J^{*} \subseteq J^{*\m} \subseteq J^{*\m^2}
\subseteq \cdots \subseteq J^{*\m^{d-1}} \subseteq J^{*\m^d} = R,
$$
where the equality on the right follows from the tight closure
Brian\c con--Skoda (\cite{HH1, Theorem 5.4}).
So it is natural to ask the following question.

\definition{Question 2.14}
Let $R$ be an excellent equidimensional local ring of characteristic $p >0$,
and let $J$ be a minimal reduction of $\m$ (in general, a parameter ideal of
$R$).  When does equality $J^{*\m^{d-1}} = J$ hold?
\enddefinition

In the following, we give a characterization of regular local rings in
terms of $\a$-tight closure, which gives an answer to the above question.
See also Section 5 about related problems.

\proclaim {Theorem 2.15}
Let $(R,\,\m,\,k)$ be an excellent equidimensional reduced local ring of
characteristic $p>0$ with $d = \dim R \ge 1$, and assume that $k$ is infinite.
Then the following conditions are equivalent.
\roster
\item $R$ is regular.
\item $\tau (\m^{d-1}) = R$, i.e., 
$I^{*\m^{d-1}}=I$ holds for every ideal $I$ of $R$.
\item $J^{*\m^{d-1}} = J$ holds for some parameter ideal $J$  of $R$.
\endroster
\endproclaim

In order to prove Theorem 2.15, we need the following lemma. Note that we
do not need to assume that $R$ is excellent in the proof of this lemma.

\proclaim {Lemma 2.16}
Let $(R,\m)$ be a regular local ring with $d = \dim R \ge 1$.
Then $I^{*\m^{d-1}} = I$ for every ideal $I$ of $R$.
\endproclaim

\demo{Proof}
Suppose that $I^{*\m^{d-1}} \ne I$ for some ideal $I$ of $R$.
Let $z \in I^{*\m^{d-1}} \setminus I$. By definition, there exists
$c \in R^{\circ}$ such that $cz^q\m^{(d-1)q} \subseteq I^{[q]}$
for all $q = p^e,\,e \gg 0$. Since the Frobenius map $F: R\to R$ is
flat by Kunz' theorem (\cite{Ku}), we have
$$
c\m^{(d-1)q} \subseteq I^{[q]} : z^q = (I:z)^{[q]} \subseteq \m^{[q]}
$$
and hence
$$c \in \m^{[q]} : \m^{(d-1)q} = \m^{q-d+1}
$$
for all $q = p^e,\,e \gg 0$ by Lemma 2.11. This is a contradiction. 
\enddemoo

\demo{Proof of Theorem $2.15$}
Note that $R$ is approximately Gorenstein by our assumption.
$(1) \Rightarrow (2)$ follows from Lemma 2.16.
Also, $(2) \Rightarrow(3)$ is trivial.

Let us prove $(3) \Rightarrow (1)$.
Take a parameter ideal $J$ such that $J^{*\m^{d-1}} = J$.
Since $J \subseteq J^{*} \subseteq J^{*\m^{d-1}}$, we have $J=J^{*}$.
Hence $R$ is F-rational and thus is Cohen--Macaulay.
Then for a minimal reduction $\q$ of $\m$, we have $\q^{*\m^{d-1}} = \q$
(see 1.14). Thus we may assume that $J$ is a minimal reduction of $\m$.
Then by virtue of Corollary 2.8, we have $\m \subseteq J$.
This implies that $R$ is regular, as required.
\enddemoo

\head 3.
Interpretation of multiplier ideals via $\a$-tight closure \endhead

We posed Theorem 2.1 as a prime characteristic analogue of Lipman's "modified
Brian{\c c}on--Skoda theorem" \cite{Li}. The original form of "modified
Brian\c con--Skoda" is stated in terms of what is called "adjoint ideals" by
Lipman. Recently, this notion is reformulated in the theory of "multiplier
ideals" from a different point of view and plays an important role in
birational algebraic geometry; see \cite{Ei}, \cite{La}. Actually, one can
define multiplier ideals with "rational coefficients"; cf.~ Section 6.

\definition{Definition 3.1}
Let $Y$ be a normal $\Q$-Gorenstein variety over a field of characteristic
zero and let $\a \subset \OO_Y$ be a nonzero ideal sheaf. Let $f \colon X
\to Y$ be a log resolution of the ideal $\a$, that is, a resolution of
singularities of $Y$ such that the ideal sheaf $\a\OO_X$ is invertible,
say, $\a\OO_X = \OO_X(-Z)$ for an effective divisor $Z$ on $X$, and that
the union $\text{Exc}(f) \cup \text{Supp}(Z)$ of the $f$-exceptional locus
and the support of $Z$ is a simple normal crossing divisor. Given a rational
number $t \ge 0$, the {\it multiplier ideal} $\J(t\cdot\a) = \J(\a^t)$
associated to $t$ and $\a$ is defined to be the ideal sheaf
$$\J (t\cdot\a) = f_*\OO_X(\rup{K_{X/Y}-tZ})$$
in $\OO_Y$, where the $\Q$-divisor $K_{X/Y} = K_X - f^*K_Y$ is the
discrepancy of $f$. For $t = 1$, we just denote $\J(\a) := \J(1\cdot\a)$.
This definition is independent of the choice of a log resolution
$f \colon X \to Y$ of $\a$.
\enddefinition

\demo{Remark $3.2$}
(1) If $\bb = \a^n$ for a nonnegative integer $n$, then $\J(t\cdot\bb)
= \J(tn\cdot\a)$, and this justifies the notation $\J(\a^t)$ with "formal
exponent" $t$. Henceforth we prefer the exponential notation $\J(\a^t)$
rather than $\J(t\cdot\a)$; cf.\ Section 6.

(2) Multiplier ideals $\J(\a)$ have properties similar to those of the ideals
$\tau(\a)$; see Proposition 1.11. Namely, if $\bb \subseteq \a$ then
$\J(\bb) \subseteq \J(\a)$, and if $\bb$ is a reduction of $\a$ then
the equality $\J(\bb)=\J(\a)$ holds; $\J(\a)\bb \subseteq \J(\a\bb)$
for any $\a,\bb$, and if $\bb$ is locally principal then the equality
$\J(\a)\bb = \J(\a\bb)$ holds; $\J(\a) \ne 0$ as long as $\a \ne 0$;
and if $Y$ has only log terminal singularities, then $\a \subseteq \J
(\a)$. Also, using vanishing theorems in characteristic zero, one can
prove a "modified Brian\c con--Skoda theorem" (\cite{Li, Theorem 1.4.1},
\cite{La}): If $\a$ is generated by $r$ elements, then $\J(\a^{n+r-1})
\subseteq \a^n$ for every $n \ge 0$; cf.\ Theorem 2.1. Later in Section 4,
we study more about similarity of the ideals $\tau(\a)$ and $\J(\a)$.
\enddemo

\demo{$3.3$. Reduction to prime characteristic}
It is proved in \cite{Ha2} and \cite{Sm2} that the multiplier ideal $\J(R)$
of the unit ideal in a normal $\Q$-Gorenstein ring $R$ essentially of finite
type over a field of characteristic zero coincides, after reduction to
characteristic $p \gg 0$, with the test ideal $\tau(R)$. We generalize
this result in Theorem 3.4 below.
To state the result, we have to begin with a ring $R$ and an ideal $\a$ in
characteristic zero, and reduce them to characteristic $p \gg 0$ together
with a log resolution $f \colon X \to \Spec R$ of $\a$.

Let $R$ be an algebra essentially of finite type over a field $k$ of
characteristic zero and let $\a \subseteq R$ be an ideal. One can choose
a finitely generated $\Z$-algebra $A$ contained in $k$ and a subalgebra
$R_A$ of $R$ essentially of finite type over $A$ such that the natural map
$R_A\otimes_Ak \to R$ is an isomorphism and $\a_A = \a \cap R_A$ generates
the ideal $\a$ of $R$. For a maximal ideal $\mu$ of $A$, we consider the
base change to its residue field $\k = \k(\mu)$ over $A$ to get a prime
characteristic ring $R_{\k} = R_A \otimes_A \k$ and an ideal $\a_{\k} =
\a_AR_{\k}$. The data consisting of $\k = \k(\mu)$, $R_{\k}$, $\a_{\k}$
is considered to be a "prime characteristic model" of the original data
in characteristic zero, and we refer to such $(\k,R_{\k},\a_{\k})$ for
maximal ideals $\mu$ in a suitable dense open subset of $\Spec A$ as
"reduction to characteristic $p \gg 0$" of $(k,R,\a)$. Furthermore,
given a morphism of schemes essentially of finite type over $k$ (and
even more, a commutative diagram consisting of a finite collection
thereof), e.g., a log resolution $f \colon X \to \Spec R$ of $\a$,
we can reduce this entire setup to characteristic $p \gg 0$. (See
\cite{Ha1}, \cite{Ha2}, \cite{HH3}, \cite{Sm1}, \cite{Sm2} for more
details.) We use the phrase "in characteristic $p \gg 0$" when we speak
of such a setup reduced from characteristic zero to characteristic
$p \gg 0$.

The main result of this section is the following theorem, which ensures the
correspondence of the ideal $\tau(\a)$ and the multiplier ideal $\J (\a)$.
See Theorem 6.7 for a generalization of this theorem to the case of
"rational coefficients."
\enddemo

\proclaim{Theorem 3.4}
Let $R$ be a normal $\Q$-Gorenstein local ring essentially of finite type
over a field and let $\a$ be a nonzero ideal. Assume that $\a \subseteq R$
is reduced from characteristic zero to characteristic $p \gg 0$, together
with a log resolution of singularities $f \colon X \to Y = \Spec R$ such
that $\a\OO_X = \OO_X(-Z)$ is invertible. Then
$$\tau(\a) = H^0(X,\OO_X(\rup{K_{X/Y}}-Z)).$$
\endproclaim

\definition{3.5}
The remainder of this section is devoted to proving Theorem 3.4.
Our strategy is to reduce to the case where the ring $R$ is
quasi-Gorenstein by passing to a canonical covering; see \cite{T} for
a direct proof which does not use a canonical covering. Let $(R,\m)$ be
a normal $\Q$-Gorenstein local ring with a canonical module $\omega_R$,
and let $r$ be the least positive integer such that the $r$th symbolic
power $\omega_R^{(r)}$ of $\omega_R$ is isomorphic to $R$. Given a fixed
isomorphism $\omega_R^{(r)} \cong R$, one has a natural ring structure
of $S = \bigoplus_{i=0}^{r-1}\omega_R^{(i)}$. This is a quasi-Gorenstein
local ring (i.e., $\omega_S \cong S$) with the maximal ideal $\m \oplus
\bigoplus_{i=1}^{r-1} \omega_R^{(i)}$, and we call $S$ a {\it canonical
covering} of $R$.

The following two lemmas make it possible to reduce the proof of the
theorem to the quasi-Gorenstein case.
\enddefinition

\proclaim{Lemma 3.6}
Let $(R,\m)$ be a normal local ring of characteristic $p > 0$, $S$ a
canonical covering of $R$ as above, and assume that $r$ is not divisible
by $p$. Then for any ideal $\a \subseteq R$,
$$\tau(\a S) \cap R = \tau (\a).$$
\endproclaim

\demo{Proof} Identical to the case where $\a = R$; see \cite{Ha2,
Section 2}, and also \cite{Sm2}.
\enddemoo

The following lemma is also proved entirely as in the same way as the
case where $\a = R$ \cite{Sm2, Proposition 3.2}, but we include the proof
for the sake of completeness.

\proclaim{Lemma 3.7 {\rm cf.\ (\cite{Sm})}}
Let $(R,\m)$ be a normal local ring essentially of finite type over a field
of characteristic zero and let $S$ be a canonical covering of $R$ as above.
Then for any ideal $\a \subseteq R$ and any rational number $t \ge 0$,
$$\J ((\a S)^t) \cap R = \J (\a^t).$$
\endproclaim

\demo{Proof}
Let $\pi \colon \Spec S \to \Spec R$ be the canonical covering, and let
$f \colon X \to \Spec R$ and $g \colon Y \to \Spec S$ be log resolutions
of $\a$ and $\a S$, respectively, which make the following diagram commute.
$$
\CD
         Y         @>g>>  \Spec S  \\
@V{\tilde{\pi}}VV       @V{\pi}VV \\
         X         @>f>>  \Spec R
\endCD
$$
Let $E$ be the reduced divisor on $X$ supported on $\text{Exc}(f)
\cup \text{Supp}(Z)$ and let $G$ be the reduced divisor with the same
support as $\tilde{\pi}^*E$. Then the log ramification formula (see e.g.,
\cite{Ka}) tells us that
$$K_{Y/S} + G = \tilde{\pi}^*(K_{X/R} + E) + P \tag{3.7.1}$$
for some effective divisor $P$ on $Y$ such that codim$(\tilde{\pi}(P),X)
\ge 2$.

Now let $u \in \J(\a^t) = H^0(X,\OO_X(\rup{K_{X/R}-tZ}))$, i.e.,
$\Div_X(u)+\rup{K_{X/R}-tZ} \ge 0$. Since $\Supp (K_{X/R}-tZ) \subseteq
E$, one has $K_{X/R} - tZ + E \ge \rup{K_{X/R} - tZ}$, so that
$$\Div_X(u) + K_{X/R} + E - tZ \ge 0,$$
and this is a strict inequality for the coefficient in each irreducible
component of $E$. Pulling this back by $\tilde\pi$ and applying (3.7.1)
give an inequality
$$\Div_Y(u) + K_{Y/S} + G - t\tilde{\pi}^*Z \ge 0,$$
which is a strict inequality for the coefficient in each irreducible
component of $G$. Since $G$ is reduced, it follows that
$\Div_Y(u) + \rup{K_{Y/S} - t\tilde{\pi}^*Z} \ge 0$.
Hence $u \in H^0(X,\OO_Y(\rup{K_{Y/S}-t\tilde{\pi}^*Z})) = \J((\a S)^t)$.

Conversely, let $u \in \J((\a S)^t) \cap R$ and fix any prime divisor
$D$ on $X$. To prove $u \in \J(\a^t)$, it is enough to show that
$a := v_D(u) + \text{coeff}_D(\rup{K_{X/R} - tZ})$ is non-negative.
If $D$ is not $f$-exceptional, this follows because $\tilde\pi$ is \'etale
and finite at the generic point of $D$ and
$\Div_Y(u)+\rup{K_{Y/S}-t\tilde{\pi}^*Z} \ge 0$ by $u \in \J((\a S)^t)$.
Now let $D$ be $f$-exceptional, $F$ any prime divisor on $Y$ dominating
$D$ (this in particular implies that $F \not\subseteq \Supp (P)$) and let
$e := \text{coeff}_F(\tilde{\pi}^*D)$. Then
$$v_F(u) + \text{coeff}_F(\tilde{\pi}^*(\rup{K_{X/R}-tZ})) = ea.
                                                         \tag{3.7.2}$$
On the other hand, since $\Div_Y(u)+\rup{K_{Y/S}-t\tilde{\pi}^*Z} \ge 0$,
it follows from (3.7.1) that
$\Div_Y(u) + \rup{\tilde{\pi}^*(K_{X/R}+E)-G+P-t\tilde{\pi}^*Z} \ge 0$.
Since $F \not\subseteq \Supp (P)$, this implies that
$$v_F(u) + \text{coeff}_F(\rup{\tilde{\pi}^*(K_{X/R} - tZ)})
        \ge \text{coeff}_F(-\tilde{\pi}^*E + G) = -e + 1. \tag{3.7.3}$$
It follows from (3.7.2) and (3.7.3) that $ea \ge -e+1$, so that $a \ge 0$,
as required.
\enddemoo

\demo{Remark}
The proof of Lemma 3.7 does not only work for canonical coverings but
also does under a weaker assumption that $R \hookrightarrow S$ is finite
and \'etale in codimension 1. It is desirable that Lemma 3.6 also holds
for finite extensions which is \'etale in codimension 1, and this issue
has been recently settled by Takagi (\cite{HT}, \cite{T}).
\enddemo

\proclaim{Proposition 3.8}
Let $(R,\m)$ be a $d$-dimensional normal local ring of characteristic $p > 0$
and let $\a$ be a nonzero ideal. Let $f \colon X \to \Spec R$ be a proper
birational morphism from a normal scheme $X$ such that $\a\OO_X = \OO_X(-Z)$
is an invertible sheaf, and denote the closed fiber of $f$ by $E= f^{-1}(\m)$.
Then one has an inclusion
$$\Ker\left( H^d_{\m}(R) \overset{\delta}\to\lra H^d_E(\OO_X(Z))\right)
\subseteq 0^{*\a}_{H_{\m}^d(R)},$$
where $\delta \colon H^d_{\m}(R) \to H^d_E(\OO_X(Z))$ is an edge map
$H^d_{\m}(R) \to H^d_E(\OO_X)$ of the spectral sequence
$H^i_{\m}(R^jf_*\OO_X) \Longrightarrow H^{i+j}_E(\OO_X)$ followed by the
natural map $H_E^d(\OO_X) \to H_E^d(\OO_X(Z))$.
\endproclaim

\demo{Proof}
First note that, for any $q = p^e$ and any $c \in \a^q \subseteq
H^0(X,\OO_X(-qZ))$, we have the following commutative diagram with
exact rows.
$$\CD
0 @>>> \Ker(\delta) @>>> H_{\m}^d(R) @>>> H_E^d(\OO_X(Z)) @>>> 0 \\
   & &      @VVV          @V{cF^e}VV        @V{cF^e}VV &          \\
0 @>>> \Ker(\delta) @>>> H_{\m}^d(R) @>>> H_E^d(\OO_X(Z)) @>>> 0
\endCD
$$
Then $\a^qF^e(\Ker(\delta)) \subseteq \Ker(\delta)$ for all $q = p^e$,
and the conclusion follows from Proposition 1.15.
\enddemoo

By virtue of Lemmas 3.6 and 3.7, it is sufficient to prove Theorem 3.4 in
the case where $R$ is quasi-Gorenstein, i.e., $\omega_R \cong R$. In this
case, however, the assertion of Theorem 3.4 coincides with the following

\proclaim{Theorem 3.9}
Let $(R,\m)$ be a $d$-dimensional normal local ring essentially of finite
type over a field of characteristic $p$ and let $\a$ be a nonzero ideal.
Assume that $\a \subseteq R$ is reduced from characteristic zero to
characteristic $p \gg 0$, together with a resolution of singularities
$f \colon X \to Y = \Spec R$ such that $\a\OO_X = \OO_X(-Z)$ is invertible.
Then
$$0^{*\a}_{H_{\m}^d(R)} =
\Ker\left( H^d_{\m}(R) \overset{\delta}\to\lra H^d_E(\OO_X(Z))\right),$$
where $E$ is the closed fiber of $f$ and $\delta$ is the edge map as
in $3.8$, or dually,
$$\Ann_{\omega_R}(0^{*\a}_{H_{\m}^d(R)}) = H^0(X,\omega_X(-Z))
                                               \quad \text{in }\omega_R.$$
\endproclaim

\demo{Proof}
First, let us discuss the situation in characteristic zero before reduction
to characteristic $p \gg 0$. In characteristic zero, we choose a nonzero
element $c \in \a$ such that $R_c$ is regular and a log resolution $f \colon
X \to \Spec R$ of the ideal $c\a$. Then $\a\OO_X = \OO_X(-Z)$ for an effective
divisor $Z$ on $X$ and $\Supp (Z+\Div_X(c))$ is a simple normal crossing
divisor. We choose an $f$-ample $f$-exceptional $\Q$-Cariter divisor $D$ and
a sufficiently small rational number $\varepsilon >0$ such that $\rd{\tilde{Z}
+ \varepsilon\Div_X(c)} = Z$, where $\tilde{Z} = Z - D$.

Now we reduce the entire setup as above to characteristic $p \gg 0$,
and switch the notation to denote things after reduction modulo $p$.
Let the ideal $\a$ be generated by $r$ elements. Since $R_c$ is regular,
some power $c^s$ of $c \; (\in \a)$
is a usual test element and also an $\a$-test element. Also, since
$-\tilde{Z}$ is $f$-ample,
$\Cal{K} = \bigoplus_{n>0} H^0(X,\omega_X(-\rd{n\tilde{Z}}))$ is a
finitely generated module over
$\Cal{R} = \bigoplus_{n\ge 0} H^0(X,\OO_X(-n\tilde{Z}))$.
Say $\Cal K$ is generated in degree $\le n_0$.

Since we are working in characteristic $p \gg 0$, the $e$-times iterated
Frobenius map
$$F^e \colon H_E^d(\OO_X(\tilde{Z})) \to
        H_E^d(\OO_X(p^e(\tilde{Z}+\varepsilon\Div_X(c)))$$
is injective for all $q = p^e$ by Proposition 3.6 and Corollary 3.8 of
\cite{Ha1}; see also \cite{Ha2, Discussion 4.6} and \cite{MS}. This implies
that the map
$$c^mF^e \colon H_E^d(\OO_X(Z))
        = H_E^d(\OO_X(\tilde{Z})) \to H_E^d(\OO_X(p^e\tilde{Z}))$$
is injective for all sufficiently large $e \in \N$ such that $p^e
\varepsilon \ge m := r + 2s + n_0 - 1$. For such $q = p^e \gg 0$, we
consider the following commutative diagram with exact rows.
$$\CD
0 @>>>  \Ker (\delta) @>>>  H_{\m}^d(R)
                                                @>{\delta}>> H_E^d(\OO_X(Z)) @>>> 0         \\
   & &      @VVV           @V{c^mF^e}VV @V{c^mF^e}VV &               \\
0 @>>> \Ker (\delta_e) @>>> H_{\m}^d(R)
                                                @>{\delta_e}>> H_E^d(\OO_X(q\tilde{Z})) @>>> 0
\endCD
$$
Here, $\Ker(\delta_e)$ is considered to be the annihilator in
$H_{\m}^d(R)$ of $H^0(X,\omega_X(-\rd{q\tilde{Z}}))$ viewed as a
submodule of $\omega_R$, with respect to the duality pairing $\omega_R
\times H_{\m}^d(R) \to H_{\m}^d(\omega_R) \cong E_R(R/\m)$.

Now, if $\xi \in H_{\m}^d(R)$ is not in $\Ker (\delta)$, then
$$c^m\xi^q \notin \Ker (\delta_e)
            = \Ann_{H_{\m}^d(R)} H^0(X,\omega_X(-\rd{q\tilde{Z}}))$$
for all $q = p^e \gg 0$ by the above commutative diagram.
Then for all sufficiently large $q = p^e (\ge n_0)$, there exists an
integer $n$ with $0 \le n \le n_0$ such that
$$c^m\xi^q \notin [0:H^0(X,\OO_X((n-q)\tilde{Z}))]_{H_{\m}^d(R)},$$
since $\Cal{K}$ is generated in degree $\le n_0$ as a graded $\Cal
R$-module. Hence it follows from $H^0(X,\OO_X((n-q)\tilde{Z}))
\subseteq H^0(X,\OO_X((n-q)Z)) \subseteq \overline{\a^{q-n_0}}$ that
$c^m\xi^q\overline{\a^{q-n_0}} \ne 0$.
On the other hand, since $c \in \a$ and $c^s \in \tau(R)$,
$$c^{m-s}\overline{\a^{q-n_0}} = c^{r+s+n_0-1}\overline{\a^{q-n_0}}
        \subseteq \tau(R)\overline{\a^{q+r-1}} \subseteq \a^{q}$$
by tight closure Brian\c con--Skoda. Thus we have that
$c^s\xi^q\a^{q} \ne 0$. But this implies that
$\xi \notin 0^{*\a}_{H^d_{\m}(R)}$,
since $c^s$ is an $\a$-test element (cf.\ Theorem 1.13).

Consequently, we have $0^{*\a}_{H^d_{\m}(R)} \subseteq \Ker (\delta)$.
The reverse inclusion follows from Proposition 3.8, and we are done.
\enddemoo


\head 4.
Properties of the ideal $\tau(\a)$ analogous to multiplier ideals \endhead

In this section we prove some properties of the ideal $\tau(\a)$
analogous to those of the multiplier ideal $\J(\a)$, which are found
in \cite{DEL}, \cite{How}, and in Lazarsfeld's lecture note \cite{La};
see also similar results for "tight closure for pairs" in \cite{T}.
The results in this section are proved in any fixed characteristic $p >
0$. However, in view of Theorem 3.4, we can also say that they provide
characteristic $p$ proofs of the properties of multiplier ideals in
characteristic zero.

We also want to remark that the results in this section can be generalized
to "rational coefficients" (see Section 6), as well as those for multiplier
ideals.

\proclaim{Theorem 4.1 {\rm (Restriction theorem, cf.\ \cite{La})}}
Let $(R,\m)$ be a normal $\Q$-Gorenstein complete local ring 
of characteristic $p > 0$ and let $x \in \m$ be a non-zerodivisor of $R$.
Let $S = R/xR$ and assume that $S$ is normal. Then for any ideal $\a$ of $R$,
$$\tau(\a S) \subseteq \tau(\a)S.$$
\endproclaim

\demo{Proof}
Let $E_R = E_R(R/\m)$ and $E_S = E_S(S/\m S)$ be the injective envelopes
of residue fields of $R$ and $S$, respectively. Then one has $E_S \cong
(0:x)_{E_R} \subset E_R$. We first prove the following claim, viewing
$E_S$ as a submodule of $E_R$ via this inclusion.

\proclaim{Claim 4.1.1}\quad
$0^{*\a}_{E_R} \cap E_S \subseteq 0^{*\a S}_{E_S}$.
\endproclaim

\demo{Proof of Claim $4.1.1$}
Since $R$ is normal, we can choose an $\a$-test element $c$ whose image
$\bar c$ in $S$ is nonzero by Proposition 1.7. 
Then we have the following commutative diagram for each $q = p^e$ (see
the proof of \cite{HW, Theorem 4.9}).
$$\CD
             E_S             @>\alpha>>          E_R        \\
     @VV \bar{c}F_S^e V                   @VV cx^{q-1}F_R^e V \\
         \F_S^e(E_S)         @>\beta >>      \F_R^e(E_R)
\endCD
$$
The map $\alpha \colon E_S \to E_R$ is the inclusion map mentioned above.
Note also that $\F_R^e(E_R) \cong H^d_{\m}(\omega_R^{(q)})$ and $\F_S^e(E_S)
\cong H^{d-1}_{\m}(\omega_S^{(q)})$ by \cite{Wa}, and $\beta \colon \F_S^e(E_S)
\to \F_R^e(E_R)$ arises as a connecting homomorphism of the long exact
sequence
$$\cdots \to H_{\m}^{d-1}(\omega_R^{(q)}) @>x>> H_{\m}^{d-1}(\omega_R^{(q)})
          \to H_{\m}^{d-1}(\omega_S^{(q)}) @>\beta>> H_{\m}^d(\omega_R^{(q)})
          @>x>> H_{\m}^d(\omega_R^{(q)}) \to 0$$
associated to $0 \to \omega_R^{(q)} @>x>> \omega_R^{(q)} \to \omega_R^{(q)}/
x\omega_R^{(q)} \to 0$ (\cite{HW}). It then follows that $\Ker(\beta) \cong
H^{d-1}_{\m}(\omega_R^{(q)})/xH^{d-1}_{\m}(\omega_R^{(q)})$ is a proper
$S$-submodule of $H_{\m}^{d-1}(\omega_S^{(q)})$. Hence we can show as in the
proof of Proposition 1.15 that $\Ker(\beta)$ is annihilated by an element
$\bar{d} \in S^{\circ}$.

Now let $\xi \in 0^{*\a}_{E_R} \cap E_S$. Then $cF_R^e(\xi)\a^q = 0$ in
$\F_R^e(E_R)$ for all $q = p^e$, since $c$ is an $\a$-test element. This
implies that $\bar{c}F_S^e(\xi)\a^q \subseteq \Ker(\beta)$ by the above
commutative diagram. Therefore $\bar{c}\bar{d}F_S^e(\xi)(\a S)^q = 0$ for
all $q = p^e$ with $\bar{c}\bar{d} \in S^{\circ}$, whence $\xi \in 0^{*\a
S}_{E_S}$, as claimed.
\enddemoo

We continue the proof of Theorem 4.1. Since $R$ is complete, $0^{*\a}_{E_R}
= (0:\tau(\a))_{E_R}$, so that
$$0^{*\a}_{E_R} \cap E_S = (0:\tau(\a)+xR)_{E_R}
                          = \left( 0:\frac{\tau(\a)+xR}{xR} \right)_{E_S}
                          = (0:\tau(\a)S)_{E_S}.$$
Hence by Claim 4.1.1, we conclude
$\tau(\a S) \subseteq \Ann_S(0^{*\a}_{E_R} \cap E_S) = \tau(\a)S$.
\enddemoo

\demo{Remark $4.2$}
Theorem 4.1 implies that F-regularity "deforms" in $\Q$-Gorenstein rings
\cite{AKM}. Namely, If $R$ is $\Q$-Gorenstein, $x \in \m$ is a nonzero-divisor
and if $S = R/xR$ is F-regular, then $R$ is also F-regular. This result fails
in the absence of $\Q$-Gorensteinness \cite{Si}. On the other hand, the
F-regularity of $R$ does not imply the F-regularity of $S = R/xR$, and
this suggests that the containment $\tau(\a S) \subseteq \tau(\a)S$ is
far from the equality in general; see also \cite{HW, Theorem 4.9}.
\enddemo

Our next objective is to apply Theorem 4.1 to showing the property of
$\tau(\a)$ called "subadditivity" (Theorem 4.5), which is established
for multiplier ideals in regular rings by Demailly, Ein and Lazarsfeld
\cite{DEL}. Our strategy is to mimic the idea of "restriction to the
diagonal" used in \cite{DEL}. We first prove the following

\proclaim{Lemma 4.3}
Let $R = k[[x_1, \dots ,x_r]]$ and $S = k[[y_1, \dots ,y_s]]$ be complete
regular local rings with residue field $k$ and let $T = R \hat{\otimes}_k S
= k[[x_1, \dots ,x_r,y_1, \dots ,y_s]]$ be their complete tensor product.
Then for any ideals $\a \subset R$ and $\bb \subset S$, we have that
$(\a \otimes S) \cap (R \otimes \bb) = (\a \otimes \bb)T$ in $T$.
\endproclaim

\demo{Proof}
We regard $R$ and $S$ as subrings of $T$ via the natural ring homomorphisms
$R \hookrightarrow T$ and $S \hookrightarrow T$. Then what is to be proved is
that $\a T \cap \bb T = \a \bb T$. To prove this we may assume without loss
of generality that $\bb$ is a proper ideal of $S$.

First we note that the composition $R \to T \to T/\bb T$ of the natural maps
is a flat ring extension. Indeed, $R = R \otimes_k k \to R \otimes_k S/\bb$
is flat since so is $k \to S/\bb$ by $\bb \cap k = 0$, and the completion
map $R \otimes_k S/\bb \to T/\bb T$ is also flat.

Now let $F_{\bullet}$ be an $R$-free resolution of $R/\a$. Then $F_{\bullet}
\otimes_R T$ is a $T$-free resolution of $T/\a T$ since $T$ is flat over $R$.
Since $T/\bb T$ is also flat over $R$, one has that
$\text{Tor}^T_i(T/\a T,T/\bb T) = H_i((F_{\bullet}\otimes_R T)\otimes_T T/\bb)
= H_i(F_{\bullet}\otimes_R T/\bb T) = 0$ for $i > 0$.
In particular,
$(\a T \cap \bb T)/\a\bb T \cong \text{Tor}^T_1(T/\a T,T/\bb T) = 0$.
Thus we conclude that $\a T \cap \bb T = \a\bb T$, as required.
\enddemoo

\proclaim{Proposition 4.4}
Let $k$ be a field of characteristic $p$ and let $R = k[[x_1, \dots ,x_r]]$,
$S = k[[y_1, \dots ,y_s]]$ and $T = R \hat{\otimes}_k S = k[[x_1, \dots ,x_r,
y_1, \dots ,y_s]]$ be as in Lemma $4.3$. Then for any ideals $\a \subset R$
and $\bb \subset S$,
$$\tau((\a \otimes \bb)T) \subseteq (\tau(\a) \otimes \tau(\bb))T.$$
\endproclaim

\demo{Proof}
Let us denote the injective envelopes of the residue fields of $R$, $S$, $T$
by $E_R$, $E_S$, $E_T$, respectively. Then we can describe them in terms of
inverse polynomials as $E_R = (x_1\cdots x_r)^{-1}k[x_1^{-1},\dots ,x_r^{-1}]$,
$E_S = (y_1 \cdots y_s)^{-1}k[y_1^{-1}, \dots ,y_s^{-1}]$, $E_T = (x_1\cdots
x_r y_1\cdots y_s)^{-1}k[x_1^{-1},\dots ,x_r^{-1},y_1^{-1},\dots ,y_s^{-1}]$,
so in particular, $E_T = E_R \otimes_k E_S$. Then it is easy to see that
$0^{*\a}_{E_R} \otimes E_S + E_R \otimes 0^{*\bb}_{E_S}
\subseteq 0^{*(\a\otimes\bb)T}_{E_T}$. 
Hence
$$
\align
\tau((\a \otimes \bb)T) & \subseteq
   \Ann_T(0^{*\a}_{E_R} \otimes E_S) \cap \Ann_T(E_R \otimes 0^{*\bb}_{E_S}) \\
   & = (\tau(\a) \otimes S) \cap (R \otimes \tau(\bb))
     = (\tau(\a) \otimes \tau(\bb))T
\endalign
$$
by Lemma 4.3.
\enddemoo

\proclaim{Theorem 4.5 {\rm (Subadditivity, cf.\ \cite{DEL})}}
Let $(R,\m)$ be a complete regular local ring of characteristic $p > 0$.
Then for any two ideals $\a$, $\bb$ of $R$,
$$\tau(\a\bb) \subseteq \tau(\a)\tau(\bb).$$
\endproclaim

\demo{Proof}
Let $T = R \hat{\otimes}_k R$ and let $\Delta \colon T \to R$ be the ring
homomorphism sending $x \otimes y \in T$ to $xy \in R$. If we restrict the
containment $\tau(\a\bb T) \subseteq \tau(\a)\tau(\bb)T$ in Proposition 4.4
by the diagonal map $\Delta \colon T \to R$, we immediately obtain
$$\tau(\a\bb) \subseteq \tau(\a\bb T)R \subseteq \tau(\a)\tau(\bb)$$
by virtue of Theorem 4.1.
\enddemoo

\demo{Remark $4.6$}
If $(R,\m)$ is not F-regular, then it is clear that the subadditivity
breaks down for $\a = \bb = R$. We do not know whether or not the
subadditivity always holds in F-regular rings in general.
\enddemo

\definition{4.7. Toric case}
In \cite{How}, Howald gave a combinatorial description of the multiplier
ideal $\J (\a)$ of a monomial ideal $\a$ in a polynomial ring over a field.
We show that the ideal $\tau(\a)$ has a similar description in a more general
situation, namely, $\a$ is a toric ideal of a toric ring $R$ over a field
$k$. Note that, in this case, the multiplier ideal $\J (\a)$ can be defined
even if char $k = p > 0$, because there exists a log resolution of $\a$ in
the toric category.

Let $M = \Z^d$, $N = \Hom_{\Z}(M,\Z)$ and denote the duality pairing of
$M_{\R} = M \otimes_{\Z} \R$ with $N_{\R} = N \otimes_{\Z} \R$ by $\langle
\; ,\; \rangle \colon M_{\R} \times N_{\R} \to \R$. Let $\sigma \subset N_{\R}$
be a strongly convex rational polyhedral cone and denote $\sigma^{\vee} =
\{m \in M_{\R} \,|\; \langle m,n \rangle \ge 0$ for all $n \in \sigma\}$ as
usual. Let $R = k[\sigma^{\vee} \cap M]$ be the toric ring over a field $k$
defined by $\sigma$, that is, the subring of a polynomial ring $k[x_1,\dots
,x_d]$ generated as a $k$-algebra by monomials $x^m=x_1^{m_1}\cdots x_d^{m_d}$
with $m = (m_1, \dots ,m_d) \in \sigma^{\vee} \cap M$. Also, let $D_1,\dots
,D_s$ be the toric divisors of $\Spec R$ corresponding to the primitive
generators $n_1,\dots , n_s \in N$ of $\sigma$, respectively. A toric ideal
$\a \subseteq R$ is an ideal of $R$ generated by monomials in $x_1,\dots ,x_d$.
Let $\a \subseteq R$ be a toric ideal and let $P = P(\a) \subset M_{\R}$ be
the Newton polygon of $\a$, that is, the convex hull of $\{m \in M\,|\; x^m
\in \a \}$ in $M_{\R}$. We denote the relative interior of $P$ in $M_{\R}$
by $\Int(P)$.

Now assume that $R$ is $\Q$-Gorenstein. Then there exists $w \in M_{\R}$
such that $\langle w,n_i\rangle = 1$ for $i = 1,\dots ,s$. Indeed,
since $\omega_R^{(r)}$ is principally generated for some $r \in \N$ and
$\omega_R^{(r)}$ corresponds to the divisor $-r\sum_{i=1}^s D_i$, we can
write $\omega_R^{(r)} = x^{m_0}R$ for some $m_0 \in M$ such that $\langle
m_0,n_i \rangle = v_{D_i}(x^{m_0}) = r$. Then set $w = m_0/r \in M_{\R}$.
\enddefinition

\proclaim{Theorem 4.8}
Let $R = k[\sigma^{\vee} \cap M]$ be a $\Q$-Gorenstein toric ring over
a field of characteristic $p > 0$ and let $w \in M_{\R}$ be as above.
Then for any toric ideal $\a \subseteq R$,
$$\tau(\a) = \J (\a)$$
and it is again a toric ideal. Moreover, for $m \in M$, the following
conditions are equivalent to each other.
\roster
\item $x^m \in \tau(\a)$.
\item $m + w \in \Int (P(\a))$.
\item $x^m \in \J (\a)$.
\endroster
\endproclaim

\demo{Proof}
We prove that $\tau(\a)$ is generated by monomials $x^m$ satisfying the
condition $m + w \in \Int (P(\a))$ in (2). It is essentially proved in
\cite{How} that $\J(\a)$ has the same property.

First, to simplify our argument, we note that $1 \in R^{\circ}$ is an
$\a$-test element, because toric rings are strongly F-regular. Hence,
an element $z \in E$ of the injective envelope $E = E_R(R/\m)$ of the
residue field $R/\m =k$ is in $0^{*\a}_E$ if and only if $z^q\a^q = 0$
in $\F^e(E) = {}^e\! R \otimes_R E$ for all $q = p^e$.

Next we will compute the Frobenius map $F^e\colon E \to \F^e(E)$ explicitly.
To do this we note that $\F^e(E) \cong H_{\m}^d(\omega_R^{(q)})$ for $q=p^e$
by \cite{Wa}, and $H_{\m}^d(\omega_R^{(q)})$ is $k$-dual to $\omega_R^{(1-q)}
= \bigoplus_{\langle m,n_i \rangle \ge 1-q} k \cdot x^m$.
Therefore
$$\F^e(E) = \bigoplus_{\langle m,n_i\rangle \le q-1} k \cdot x^m
           = \bigoplus_{m \in (q-1)w - \sigma^{\vee}} k \cdot x^m,$$
and the Frobenius map $F^e \colon E \to \F^e(E)$ sends $x^m \in E$ to
$x^{mq} \in \F^e(E)$.

It is now clear that $0^{*\a}_E$ and hence $\tau(\a) = \Ann_R(0^{*\a}_E)$
is generated by monomials, because everything involving is $\Z^s$-graded.

We describe the $\a$-tight closure $0^{*\a}_E$ of zero 
in $E = \bigoplus_{u\in -\sigma^{\vee}\cap M} k \cdot x^u$. Let $u \in
-\sigma^{\vee}\cap M$. Then $x^u \in 0^{*\a}_E$ if and only if $x^{qu}\a^q
= 0$ in $\F^e(E)$, or equivalently, $(qu+qP) \cap ((q-1)w-\sigma^{\vee})
\cap M = \emptyset$, for all $q=p^e$. Dividing out by $q$, we can rephrase
this into the condition that $(u+P)\cap\Int(w-\sigma^{\vee}) = \emptyset$,
because $-w/q \in \Int(-\sigma^{\vee})$. Since this is equivalent to saying
that $\Int(u + P) \cap (w - \sigma^{\vee}) = \emptyset$, it follows that
$x^u \in 0^{*\a}_E$ if and only if $\Int(P) \cap (w - u - \sigma^{\vee})
= \emptyset$, or equivalently, if $w - u \notin \Int(P)$.

Now the equivalence of conditions (1) and (2) follows immediately, because
a monomial $x^m \in R = k[\sigma^{\vee} \cap M]$ is in $\tau(\a) = \Ann_R
(0^{*\a}_E)$ if and only if $x^{-m} \notin 0^{*\a}_E$.
\enddemoo

\definition{Example 4.9}
Let $S = k[x_1,\dots ,x_d]$ be a polynomial ring and let $R = S^{(r)}$ be
the $r$th Veronese subring of $S$. We can easily compute the ideal $\tau(\a)
= \J(\a)$ associated to a monomial ideal $\a$ of $R$ as follows; cf.\
\cite{How}.

We choose $N$ and $M$ to be the overlattice $N = \Z^d+\frac1r (1,\dots ,1)
\Z$ and the sublattice $M = \{m \in \Z^d |\; \langle m,n \rangle \in \Z\}$
of $\Z^d$, respectively, and let $\sigma$ be the first orthant in $N_{\R}
= \R^d$. Then the dual cone $\sigma^{\vee}$ is also the first orthant in
$M_{\R} = \R^d$, and the ring $R = k[\sigma^{\vee} \cap M]$ is realized as
it is as the $r$th Veronese subring of $S = k[\sigma^{\vee} \cap \Z^d] =
k[x_1,\dots ,x_d]$. In this setting, the vector $w \in M_{\R}$ (for both
$R$ and $S$) defined in 4.7 is equal to ${\bold 1} = (1,1,\dots ,1)$. Also,
for a monomial ideal $\a \subseteq R$, the Newton polygons $P(\a)$ and
$P(\a S)$ are equal to each other in $M_{\R} = \R^d$. Therefore, Theorem
4.8 tells us that a monomial $x^m$ in $R$ (resp.\ $S$) is in $\tau(\a)$
(resp.\ $\tau(\a S)$) if and only if $m + {\bold 1} \in \Int(P(\a)) =
\Int(P(\a S))$ and in particular,
$$\tau(\a) = \tau(\a S) \cap R \; ;$$
cf.\ Lemma 3.6. For example, if $\a = \m^l$ is a power of the graded maximal
ideal $\m$ of $R$, we have
$\tau(\m^l) = \m^{\rup{l-(d-1)/r}} = \m^{\rd{l+1-d/r}}$.
\enddefinition

\head 5. F-rationality of Rees algebras and the behavior of $\tau(I)$
\endhead

Throughout this section, we assume that $(R,\,\m)$ is an excellent {\it
Gorenstein} local domain of characteristic $p > 0$ and that $I$ is an
$\m$-primary ideal of $R$. Put $d = \dim R \ge 2$. Let $\Rees(I) = R[It]$
denote the Rees algebra of $I$ over $R$, and $\M = \m\Rees(I)+\Rees(I)_+$,
the unique homogeneous maximal ideal of $\Rees(I)$. We will denote by
$\Rees'(I) = R[It,t^{-1}]$ and $G(I) = \Rees'(I)/t^{-1}\Rees'(I) =
\bigoplus_{n\ge 0} I^n/I^{n+1}$ the extended Rees algebra and the
associated graded ring of $I$, respectively. Also, let $\omega_{\Rees(I)}$
denote the graded canonical module of $\Rees(I)$, and let $\pi \colon Y =
\Proj \Rees(I) \to \Spec R$ be the blowing-up with respect to $I$.

The main purpose of this section is to describe $\omega_{\Rees(I)}$ in
terms of $\tau(I^n)$ under the assumption that $\Rees(I)$ is F-rational.
Actually, we prove the following theorem.

\proclaim {Theorem 5.1}
  Let $(R,\m)$ be an excellent Gorenstein local domain of characteristic
$p > 0$ with $d = \dim R \ge 2$. Let $I$ be an $\m$-primary ideal of $R$
and $J$ its minimal reduction. Then $\tau(I) \subseteq J:I^{d-1}$.
If we assume that $\Rees(I)$ is F-rational, then $\tau(I) = J:I^{d-1}$ and
$$
  \omega_{\Rees(I)} = \bigoplus_{n \ge 1} H^0(Y,I^n\omega_Y)
                    \cong \bigoplus_{n \ge 1} \,\tau(I^n).
$$
\endproclaim

\demo{Discussion $5.2$}
The above theorem is motivated by Hyry's papers \cite{Hy1}, \cite{Hy2}.
For example, the description of $\omega_{\Rees(I)}$ in Theorem 5.1
corresponds to the following fact used in \cite{Hy1}:
Let $(R,\m)$ be a regular local ring essentially of finite type
over a field of characteristic zero, and let $I$ be an ideal of $R$.
Suppose that $\Proj \Rees(I)$ has rational singularities.
The graded canonical module of $\Rees(I)$ is then
$\omega_{\Rees(I)} = \bigoplus_{n\ge 1} \J(I^n)$.
\par
Actually, if $\Rees(I)$ is F-rational, then so is $Y = \Proj \Rees(I)$,
whence $Y$ is pseudo-rational \cite{Sm1}. Therefore, if $Y$ has a
resolution of singularities $f \colon X \to Y$, then
$H^0(X,\,I^n\omega_X) \cong H^0(Y,\,I^n\omega_Y)$ for every $n \ge 0$.
The left-hand side of this equality coincides with the multiplier ideal
$\J(I^n)$ via the isomorphism $\omega_R \cong R$ as long as $\J(I^n)$
is defined. Moreover,
$[\omega_{\Rees(I)}]_n \cong H^0(Y,\,I^n\omega_Y)$
since $\Rees(I)$ is Cohen--Macaulay; see e.g.\ \cite{HHK}.
In particular, we have $\tau(I) = \J(I)$ in this case.
\par
Consequently, Theorem 5.1 claims that the F-rationality of $\Rees(I)$
gives a sufficient condition for which $\tau(I) = \J(I)$ holds in any
fixed positive characteristic.
\enddemo

We obtain the following corollary from Theorem 5.1 and
\cite{HWY1, Corollary 3.3}.

\proclaim {Corollary 5.3}
Suppose that $(R,\m)$ is a two-dimensional rational double point.
Let $I$ be an $\m$-primary integrally closed ideal of $R$ and
$J$ its minimal reduction.
Then $\tau(I) \subseteq J:I \; (=\J(I))$.
Also,
$\Rees(I)$ is F-rational if and only if $\tau(I) = J:I$.
\endproclaim

One can easily check the following example using the method developed
in \cite{HWY1, Section 3}.

\example{Example 5.4 {\rm (cf.\ \cite{HWY1, Theorem 3.1})}}
\par
(1) Let $(R,\m)$ be a two-dimensional excellent Gorenstein F-rational
local ring (i.e., F-rational double point), and $I$ an $\m$-primary
integrally closed ideal of $R$.
Then $\tau(I) = J:I$ for any minimal reduction $J$ of $I$.
\medskip
(2) Let $R = k[[x,y,z]]/(x^2+y^3+z^5)$, where $k$ is an algebraically
closed field of characteristic $2$. Put $\m = (x,y,z)R$ and $J = (y,z)R$.
Then $R$ is a two-dimensional rational double point but not F-rational.
Also, we have
\roster
\item"(a)" $R(\m)$ is {\it not} F-rational.
\item"(b)" $J^{*\m} = \m$. In particular, $J:J^{*\m} = \m$.
\item"(c)" $(J^{[2]})^{*\m} = (y^2,\,z^2,\,xy)$.
In particular, $J^{[2]} : (J^{[2]})^{*\m} = (x,y,z^2)$.
\item"(d)" $\tau(\m) \subseteq (x,y,z^2) \subsetneq \J(\m)=\m$.
\endroster
\endexample

In the following, we will prove Theorem 5.1.

\proclaim {Lemma 5.5}
Let $(R,\m)$ be a Gorenstein local ring of any characteristic.
Also, let $I$ be an $\m$-primary ideal of $R$ and put
$G(I) = \bigoplus_{n \ge 0} I^n/I^{n+1}$,
the associated graded ring of $I$.
Assume that $[H_{\M}^d(\Rees(I)]_0 =[H_{\M}^d(\Rees(I))]_{-1}= 0$
$($e.g.\ $\Rees(I)$ is Cohen--Macaulay$)$.
Then
$R/[\omega_{\Rees(I)}]_1 \cong \left([H_{\M}^d(G(I))]_{-1}\right)^{\vee}$,
where $(\quad)^{\vee}$ denotes the Matlis dual of $R$.
\endproclaim

\demo{Proof}
Consider the following two standard exact sequences:
$$\matrix
  0 \lra & \Rees(I)_{+} & \lra \Rees(I) \lra & R & \lra 0, \\
  \phantom{.} & & & & \\
  0 \lra & \!\!\!\! \Rees(I)_{+}(1) \!\!\!\! & \!\!\!\! \lra
    \Rees(I) \lra \!\!\!\! & \!\!\!\! G(I) \!\!\!\! & \lra 0.
\endmatrix
$$
\par
 From the first exact sequence, we have
$$
  0 = [H_{\M}^{d}(\Rees(I))]_0 \to H_{\m}^d(R) \to
      [H_{\M}^{d+1}(\Rees(I)_+)]_0 \to [H_{\M}^{d+1}(\Rees(I))]_0 = 0,
$$
where the vanishing on the right follows from $a(\Rees(I)) = -1$.
Since $R$ is Gorenstein, $R = \omega_R \cong (H_{\m}^d(R))^{\vee} \cong
\left([H_{\M}^{d+1}(\Rees(I)_{+})]_0\right)^{\vee}$.
On the other hand, the second exact sequence gives
$$
  0 = [H_{\M}^{d}(\Rees(I))]_{-1} \to [H_{\M}^d(G(I))]_{-1}
  \to [H_{\M}^{d+1}(\Rees(I)_+)]_0 \to [H_{\M}^{d+1}(\Rees(I))]_{-1} \to 0.
$$
Dualizing the above sequence, we get
$$
  0
   \lra [\omega_{\Rees(I)}]_{1}
    \lra \left([H_{\M}^{d+1}(\Rees(I)_{+})]_0\right)^{\vee} \cong R
     \lra \left([H_{\M}^d(G(I))]_{-1}\right)^{\vee}
      \lra 0.
$$
This yields the required assertion.
\enddemoo

\proclaim {Proposition 5.6}
Under the same notation as in Lemma $5.5$, assume further that $\Rees(I)$
is Cohen--Macaulay. Then $[\omega_{\Rees(I)}]_1 = J:I^{d-1}$ for every
minimal reduction $J$ of $I$.
\endproclaim

\demo{Proof}
Let $x_1,\,x_2,\,\ldots,x_d$ be a minimal system of generators of $J$.
Put $G :=G(I)$. Then it is well-known that $G(I)$ is Cohen-Macaulay and
the images in $G$ of $x_1t,\,\ldots,x_dt \; (\in \Rees(I)_1)$ form a
regular sequence (\cite{GS}). Setting $x_i^* := x_it \mod I^2$ for each
$i$, we have an exact sequence
$$
  0 \lra G(-1) \overset x_1^{*} \to \lra G
    \lra G/x_1^*G \cong G(I/x_1R) \lra 0
$$
by \cite{VV}. From this sequence, we get an exact sequence
$$
  0 = H_{\M}^{d-1}(G) \lra H_{\M}^{d-1}(G/x_1^{*}G)
    \lra H_{\M}^d(G)(-1) \overset x_1^* \to \lra H_{\M}^d(G) \lra 0.
$$
Since $a(G) \le -1$ (\cite{GS}), we have
$[H_{\M}^{d-1}(G/x_1^{*}G)]_0 \cong [H_{\M}^d(G)]_{-1}$ and
$a(G/x_1^{*}G) \le -1+1 = 0$.
By repeating the above argument, we get
$$
  [H_{\M}^d(G)]_{-1} \cong [H_{\M}^0(G/x_1^{*},\,\ldots,x_d^{*})G]_{d-1}
                     \cong [H_{\M}^0(G(I/J))]_{d-1} = \dfrac{J+I^{d-1}}{J}.
$$
Also, since $R/J$ is Gorenstein, we have
$\left(\frac{J+I^{d-1}}{J}\right)^{\vee} \cong R/J:I^{d-1}$.
Combining this with the previous lemma, we get
$$
  [\omega_{\Rees(I)}]_1 = \Ann_R \left([H_{\M}^d(G)]_{-1}\right)^{\vee}
         = \Ann_R \left(\dfrac{J+I^{d-1}}{J}\right)^{\vee} = J:I^{d-1}.
$$
\enddemoo

\proclaim {Proposition 5.7}
Let $(R,\,\m)$ be an excellent Gorenstein local domain of characteristic
$p > 0$, and let $I$ be an $\m$-primary ideal of $R$. Also, let $J$ be a
minimal reduction of $I$. Then we have the following statements.
\roster
\item $\tau(I) \subseteq J:J^{*I} \subseteq J:I^{d-1}$.
\item If $\Rees(I)$ is F-rational, then $\tau(I) = J:J^{*I} = J:I^{d-1}$.
\endroster
\endproclaim

\demo{Proof}
Let $x_1,\,\ldots,x_d$ be a system of generators of $J$ and put
$J^{[l]} := (x_1^l,\,\ldots,x_d^l)$ for all integers $l \ge 1$.
\par
(1) One has $\tau(I) \subseteq J :J^{*I} \subseteq J:(J+I^{d-1}) = J:I^{d-1}$
by the definition of $\tau(I)$ and Corollary 2.8.
\par
  To see (2), we may assume that $\Rees(I)$ is Cohen--Macaulay.
Then $I^d = JI^{d-1}$ (\cite{GS}). Hence
$I^{dl-1} = J^{dl-d}I^{d-1}
           = \left(J^{[l]}J^{dl-d-l}+(x_1\cdots x_d)^{l-1}R \right) I^{d-1}$.
Thus
$J^{[l]}+I^{dl-1} = J^{[l]}+(x_1\cdots x_d)^{l-1}I^{d-1}$.
In particular, for all $l \ge 1$, we have
$J^{[l]}:(J^{[l]}+I^{dl-1}) = (J^{[l]}:(x_1\cdots x_d)^{l-1}):I^{d-1}
                             = J:I^{d-1}$.
\par
Now suppose that $\Rees(I)$ is F-rational.
Then since $(J^{[l]})^{*I} = J^{[l]}+I^{dl-1}$ by Corollary 2.9, we
have that $J^{[l]} : (J^{[l]})^{*I} = J:I^{d-1}$ for all $l \ge 1$.
Hence $\tau(I) = $  $J:J^{*I} = J:I^{d-1}$, as required.
\enddemoo

\demo{Proof of Theorem $5.1$}
  Note that $\Rees(I^n)$ is F-rational if so is $\Rees(I)$.
Actually, it is a module-finite pure subring of $\Rees(I)$.
Thus the required assertion immediately follows from Propositions
5.6 and 5.7.
\enddemoo

\demo{Proof of Corollary $5.3$}
Let $I$ be an $\m$-primary integrally closed ideal and $J$ its
minimal reduction. Then it is well-known that $I^2 =JI$ and thus
$\Rees(I)$ is Cohen--Macaulay.
\par
It is enough to show that $\tau(I) = J:I$ implies that $\Rees(I)$
is F-rational. Suppose that $\tau(I) = J:I$. Since
$\tau(I) \subseteq J^{[l]} : (J^{[l]})^{*I}
          \subseteq J^{[l]} : (J^{[l]}+I^{2l-1}) = J:I$,
in general, we have
$J^{[l]} : (J^{[l]})^{*I} = J^{[l]} : (J^{[l]}+I^{2l-1})$.
This implies that
$( J^{[l]})^{*I} = J^{[l]} + I^{2l-1}$ for all $l \ge 1$
  because $R/J^{[l]}$ is an Artinian Gorenstein local ring.
By \cite{HWY1,\,Corollary 3.3(2)}, we conclude
that $\Rees(I)$ is F-rational.
\enddemoo

In the latter half of this section, we will give some applications
of Theorem 5.1. Before stating our results, let us recall the notion
of $a$-invariant. Let $I$ be an $\m$-primary ideal of $R$ and put
$G := G(I)$ and $\M := \m \Rees(I) + \Rees(I)_{+}$.
Then the $a$-invariant $a(G)$ of $G$ is defined by
$a(G): = \max\{n \in \Z\,|\, [H_{\M}^d(G)]_n \ne 0\}$.
See \cite{GW} for details.

\proclaim {Proposition 5.8}
Let $(R,\m)$ be an excellent Gorenstein local domain of characteristic
$p > 0$. Let $I$ be an $\m$-primary ideal of $R$.
Suppose that $\Rees(I)$ is F-rational and $G:=G(I)$ is Gorenstein.
Then $\tau(I^n) = I^{n+a(G)+1}$ for all integers $n \ge 1$.
\endproclaim

\demo{Proof}
The F-rationality of $\Rees(I)$ implies that
$\tau(I^n) = [\omega_{\Rees(I)}]_n = H^0(Y,\,I^n\omega_Y)$ for all
$n \ge 1$; see Discussion 5.2. Also, since $G$ is Cohen--Macaulay,
we have
$$
\omega_G \cong \bigoplus_{n\ge 1} H^0(Y,I^{n-1}\omega_Y)/ H^0(Y,I^n\omega_Y)
  \tag{5.8.1}
$$
and $R = H^0(Y,\omega_Y) = \cdots = H^0(Y,I^{-a-1}\omega_Y)$,
where $a = a(G) \le -1$; see e.g.\ \cite{Hy2, Theorem 2.2}.
On the other hand, as $G$ is Gorenstein, we have
$$
  \omega_G \cong G(a) = \bigoplus_{n\ge -a} I^{n+a}/I^{n+a+1} \tag{5.8.2}
$$
Comparing (5.8.1) with (5.8.2), one can easily see that
$\tau(I^n) = I^{n+a+1}$ by induction on $n \ge 1$.
\enddemoo

\proclaim {Corollary 5.9}
Let $(R,\,\m)$ be an $($excellent$)$ regular local ring of characteristic
$p > 0$. Then $\tau(\m^n) = \m^{n-d+1}$ for all $n \ge 1$.
\endproclaim

\demo{Proof}
Suppose that $R$ is a regular local ring. Then $R(\m)$ is F-rational and
$G(\m) \cong k[X_1,\,\ldots,X_d]$ is Gorenstein with $a(G(\m)) = -d$.
Hence we can apply the above proposition.
\enddemoo

\remark {Remark $5.10$}
Corollary 5.9 is a generalization of the implication (1) $\Rightarrow$
(2) in Theorem 2.15. This also follows from Theorem 4.8.
\endremark

\medskip
Let $J \subseteq I$ be ideals of $R$. 
Recall that the {\it coefficient ideal} of $I$ relative to $J$,
denoted by $\a(I,J)$, is defined to be the largest ideal $\a$ of $R$
for which $I\a = J\a$.

\remark {Remark $5.11$}
In \cite{Hy2}, Hyry proved that if $R$ is a Gorenstein local ring and
$\Rees(I)$ has rational singularities then $\J(I^{d-1}) = \a(I,J)$.
In fact, a similar result follows from Theorem $5.1$ and \cite{Hy2,
Theorem 3.4}:
Let $(R,\,\m)$ be an excellent Gorenstein local domain of
characteristic $p > 0$. Let $I$ be an $\m$-primary ideal of $R$,
and $J$ its minimal reduction. If $\Rees(I)$ is F-rational, then
$\tau(I^{d-1}) = H^0(Y,I^{d-1}\omega_Y) = \a(I,J)$.
In particular, if, in addition, $I^2 =JI$, then $\tau(I^{d-1}) = J:I$.
\endremark

\medskip
In the rest of this section, we direct our attention to the ideal
$\tau(\m)$. Let $(R,\m)$ be an excellent Gorenstein F-rational local
domain of characteristic $p > 0$.
Then $\tau(\m) \supseteq \m$, that is, $\tau(\m) = \m$ or $\tau(\m) =R$.
For example, if $R$ is a regular local ring with $\dim R \ge 2$, then
$\tau(\m) = R$. More generally, we have the following proposition.

\proclaim {Proposition 5.12}
Let $(R,\m)$ be an excellent Gorenstein local domain of characteristic
$p > 0$ with $d = \dim R \ge 2$. Suppose that there exists an $\m$-primary
ideal $I$ such that $\Rees(I)$ is F-rational with $a(G(I)) \ne -1$.
Then $R$ is F-rational with $\tau(\m) = R$.
\endproclaim

\demo{Proof}
By virtue of \cite{HWY1, Corollary 2.13}, $R$ is F-rational.
\par
By Theorem 5.1, we have $\tau(I) = J:I^{d-1}$ for any minimal reduction
$J$ of $I$. Since $\Rees(I)$ is Cohen--Macaulay with $a(G(I)) \ne -1$,
we have that $I^{d-1} = JI^{d-2} \subseteq J$. Hence $\tau(I) = R$.
In particular, $\tau(\m) = R$ because $\tau(\m) \supseteq \tau(I)$.
\enddemoo

In view of the above proposition it is natural to ask the following

\example {Question 5.13}
Let $(R,\m)$ be an excellent Gorenstein F-rational local domain of
characteristic $p>0$ with $\tau(\m) = R$ and $\dim R \ge 2$.
When is $\Rees(\m)$ F-rational then?
\endexample

In case of two-dimensional local rings, $\tau(\m) = R$ implies that
$R$ is regular. Then $R(\m)$ is Gorenstein and F-rational with
$a(G(\m)) = -2$. As for three-dimensional local rings, we have
the following answer to the above question.

\proclaim {Proposition 5.14}
Let $(R,\,\m)$ be a three-dimensional excellent Gorenstein local ring
which is not regular. Then the following conditions are equivalent.
\roster
\item $\tau(\m) = R$, that is, $I^{*\m} = I$ holds for every ideal $I$
of $R$.
\item $J^{*\m} = J$ holds for some parameter ideal $J$ of $R$.
\item $\Rees'(\m)$ is F-rational and $\m^2 = J\m$ for some minimal
reduction $J$ of $\m$.
\item $\Rees(\m)$ is Gorenstein and F-rational.
\item $\tau(\m^{n}) = \m^{n-1}$ holds for all integers $n \ge 1$.
\endroster
\endproclaim

\demo{Proof}
$(1) \Rightarrow (2)$ and $(5) \Rightarrow (1)$ are trivial.
$(4) \Rightarrow (5)$ follows from Proposition 5.8 since $a(G(\m)) =-2$
(\cite{GS}).

To see $(2) \Rightarrow (3)$, we may assume that $J$ is a minimal
reduction of $\m$; see Discussion 1.14. By Corollary 2.8, we have
$\m^2 \subseteq J$, and thus $\m^2 =J\m$.
Also, $\Rees'(\m)$ is F-rational by \cite{HWY1, Corollary 4.5}.
\par
$(3) \Rightarrow (4)$:
Note that a Gorenstein local ring having minimal multiplicity is a
hypersurface with multiplicity at most $2$. Thus $R$ and $G(\m)$ are
hypersurfaces and $a(G(\m)) = 1 - \dim R = -2$. Hence $\Rees(\m)$
is Gorenstein (\cite{GS}). Also, as $\Rees'(\m)$ is Gorenstein and
F-rational, $\Rees(\m)$ is F-regular.
\enddemoo

\demo {Discussion $5.15$}
We can generalize the equivalence of $(2)$ and $(3)$ in
Proposition $5.14$ as follows; see also Theorem $2.15$.
\par
Let $(R,\m)$ be an excellent equidimensional reduced local ring of
characteristic $p > 0$. Then the following conditions are equivalent.
\roster
\item $\Rees'(\m)$ is F-rational and $\m^2 =J\m$ for some minimal
reduction $J$ of $\m$.
\item $J^{*\m^{d-2}} = J$ holds for every (or equivalently, some)
parameter ideal $J$ of $A$.
\endroster
\par
If, in addition, $R$ is Gorenstein, then the following condition is
also equivalent to the above conditions.
\roster
\item"(3)" $\tau(\m^{d-2}) = R$.
\endroster
\enddemo

\example {Example 5.16 {\rm (cf. \cite{HWY2, Proposition 3.12},
\cite{HWY1, Sect.\ 5})}}
Let $(R,\m)$ be an excellent three-dimensional Gorenstein normal local
domain of characteristic $p > 0$. If $R$ admits a non-zerodivisor
$f \in \m$ such that $R/fR$ is F-rational, then $\tau(\m) = R$.
\par
For example, let $R = k[[x,y,z,w]]/(x^2+y^a+z^b+w^c)$, where $k$ is
a field of characteristic $p >0$ and $a,\,b,\,c$ are integers with
$2 \le a \le b\le c \ll p$.
If $1/2+1/a+1/b >1$, then $\tau(\m) = R$. Otherwise, $\tau(\m) = \m$.
\endexample

\remark {Remark $5.17$}
(1) If $R$ is a three-dimensional regular local ring, then
$\tau(\m^2) = R$ (and thus $\tau(\m)= R$) and $\Rees(\m)$ is F-rational.
But $\Rees(\m)$ is not Gorenstein and $\tau(\m^n) = \m^{n-2}$ for all
$n \ge 1$.

(2) We have no examples of Gorenstein local ring $R$ for which
$\tau(\m) = R$ but $\Rees(\m)$ is not F-rational.
\endremark

\demo{Discussion $5.18$}
Let $(R,\m)$ be a complete regular local ring of characteristic $p>0$
with $d = \dim R \ge 2$, and let $I$ be an $\m$-primary ideal of $R$.
Then we expect that $\tau(I) \supsetneq I$.
\par
For example, this is true if $R(I)$ is F-rational.
We sketch a proof here. Suppose that $\tau(I) = I$.
Then $\tau(I^n) \subseteq \tau(I)^n = I^n$ for all $n \ge 1$ by the
subadditivity (Theorem 4.5). On the other hand, since $R$ is F-regular,
we have $\tau(I^n) \supseteq I^n$ in general. Also, by Theorem 5.1,
we have $\omega_{\Rees(I)} = \bigoplus_{n \ge 1} I^n = \Rees(I)_+$.
In particular, since $\Rees(I)/\omega_{\Rees(I)} \cong R$ is regular,
so is $\Rees(I)$. (Note: Recently, S.~Goto et.~al.\ proved a more
general result.) But this is impossible because $\dim R \ge 2$.
Hence $\tau(I) \supsetneq I$, as required.

As for multiplier ideals, the authors are informed of the following
result by K.-i.~Watanabe: Let $R$ be a regular local ring essentially
of finite type over a characteristic zero.
Then $\J(I) \supsetneq I$ for any $\m$-primary ideal $I$ of $R$.
\enddemo

\head 6. Rational coefficients \endhead

Recently, the theory of multiplier ideals with "rational coefficients" has
been developed and applied successfully to various problems in algebraic
geometry and commutative algebra (\cite{ELS}, \cite{La}). This motivates
us to extend the notions of $\a$-tight closure and the ideal $\tau(\a)$ to
those with "rational coefficients." In this last section we make a few
remarks on rational coefficients and address the results which generalize
in this setting.

\definition{Definition 6.1}
Let $\a$ be an ideal of a Noetherian ring $R$ of characteristic $p > 0$
such that $\a \cap R^{\circ} \ne \emptyset$ and let $N \subseteq M$ be
$R$-modules. Given a rational number $t \ge 0$, the {\it $t\cdot\a$-tight
closure} $N^{*t\cdot\a}_M$ (or, {\it $\a^t$-tight closure} $N^{*\a^t}_M$,
see Remark 6.2 (1) below) of $N$ in $M$ is defined to be the submodule of
$M$ consisting of all elements $z \in M$ for which there exists $c \in
R^{\circ}$ such that
$$cz^q\a^{\rup{tq}} \subseteq N^{[q]}_M$$
for all $q = p^e \gg 0$, where $\rup{tq}$ denotes the least integer which
is greater than or equal to $tq$.
\enddefinition

\remark{Remark $6.2$}
(1) Definition 6.1 does not change if we replace "$cz^q\a^{\rup{tq}}
\subseteq N^{[q]}_M$" (rounding up $tq$) by "$cz^q\a^{\rd{tq}} \subseteq
N^{[q]}_M$" (rounding down $tq$), as long as $\a \cap R^{\circ} \ne
\emptyset$. This is because the difference of $\rup{tq}$ and $\rd{tq}$
as the exponents of $\a$ is "absorbed" by the term $c \in R^{\circ}$.
Similarly, it is easy to see that $t\cdot\a^n$-tight closure is the
same as $tn\cdot\a$-tight closure for every nonnegative integer $n$;
cf.\ the proof of \cite{HW, Proposition 2.6}. This being so, it is
preferable to say "$\a^t$-tight closure" rather than "$t\cdot\a$-tight
closure." In the sequel, we always use "exponential notation" in this
manner and denote the $\a^t$-tight closure of $N$ in $M$ by $N^{*\a^t}_M$.

(2) We can even define tight closure with several rational coefficients
(or, several rational exponents). Namely, given ideals $\a_1,\dots ,\a_n
\subseteq R$ and rational numbers $t_1,\dots ,t_n \ge 0$, the $\a_1^{t_1}
\cdots\a_n^{t_n}$-tight closure $N^{*\a_1^{t_1}\cdots\a_n^{t_n}}_M$ is
defined by replacing "$cz^q\a^{\rup{tq}} \subseteq N^{[q]}_M$" in Definition
6.1 by "$cz^q\a_1^{\rup{t_1q}}\cdots\a_n^{\rup{t_nq}} \subseteq N^{[q]}_M$."

We also have an analogous notion of $\Delta$-tight closure for a pair
$(R,\Delta)$ of a normal ring $R$ and a $\Q$-Weil divisor $\Delta$ on
$Y = \Spec R$; see \cite{T}, \cite{HW}. If $\a_i = x_iR$ and $\Delta =
\sum_{i=1}^n t_i\cdot\Div_Y(x_i)$ for $x_i \in R$ and $0 \le t_i \in \Q$
with $1 \le i \le n$, then $\a_1^{t_1}\cdots\a_n^{t_n}$-tight closure
is the same as $\Delta$-tight closure.
\endremark

\definition{Definition 6.3}
Let $\a$ be an ideal of a Noetherian ring $R$ of characteristic $p>0$
such that $\a\cap R^{\circ} \ne \emptyset$ and let $t \ge 0$ be a
rational number. We say that an element $c \in R^{\circ}$ is an {\it
$\a^t$-test element} if $cz^q\a^{\rup{tq}} \subseteq I^{[q]}$ for all
$q = p^e$ whenever $z \in I^{*\a^t}$.
\enddefinition

Some results for $\a$-tight closure generalize to those for $\a^t$-tight
closure without essential change of proofs. However, we must be careful
about the difference of round-up and round-down when speaking of $\a^t$-test
elements. As a matter of fact, the following theorem is proved similarly
as Theorem 1.7 (1), but the proof does not work if we replace $\rup{tq}$
by $\rd{tq}$ in Definition 6.3.

\proclaim{Theorem 6.4}
Let $R$ be an F-finite reduced ring of characteristic $p > 0$ and let
$c \in R^{\circ}$ be an element such that the localized ring $R_c$ is
strongly F-regular. Then some power $c^n$ of $c$ is an $\a^t$-test
element for all ideals $\a \subseteq R$ with $\a \cap R^{\circ} \ne
\emptyset$ and all rational numbers $t \ge 0$.
\endproclaim

We can define the ideal $\tau(\a^t)$ in a similar way as in
Proposition-Definition 1.9. 
Also, Theorem 1.13 generalizes to the case of $\tau(\a^t)$ with the
same proof. We summarize the results for excellent reduced local rings
in the following.

\proclaim{Definition-Theorem 6.5}
Let $(R,\m)$ be an excellent reduced local ring of characteristic
$p > 0$ with $E = E_R(R/\m)$ and let $\a \subseteq R$ be an ideal such
that $\a\cap R^{\circ} \ne \emptyset$. Given a rational number $t \ge 0$,
we define the ideal $\tau(\a^t) \subseteq R$ by
$$\tau(\a^t) = \bigcap_M \Ann_R(0^{*\a^t}_M)
              = \bigcap_{M\subseteq E} \Ann_R(0^{*\a^t}_M)
              = \bigcap_{I\subseteq R} (I:I^{*\a^t}),$$
where $M$ runs through all finitely generated $R$-modules $($resp.\
finitely generated $R$-submodules of $E)$ in the second term $($resp.\
the third term$)$, and $I$ runs through all ideals of $R$. Moreover,
if $R$ is normal and $\Q$-Gorenstein, then
$$\tau(\a^t) = \Ann_R(0^{*\a^t}_E).$$
\endproclaim

\demo{Remark $6.6$}
We can define the ideal $\tau(\a_1^{t_1}\cdots\a_n^{t_n})$ with several
rational coefficients by replacing $\a^t$-tight closure in 6.5 by
$\a_1^{t_1}\cdots\a_n^{t_n}$-tight closure as defined in Remark 6.2 (2).
See Theorem 6.10 (2).
\enddemo

Proposition 1.15 also generalizes without essential change of the proof,
but we cannot replace the round-up $\rup{tq}$ by the round-down $\rd{tq}$
in the following

\proclaim{Proposition 6.7}
Let $(R,\m)$ be a $d$-dimensional excellent normal local ring of
characteristic $p>0$, $\a \subseteq R$ an ideal such that $\a\cap R^{\circ}
\ne \emptyset$ and let $t \ge 0$ be a rational number.
Then $0^{*\a^t}_{H^d_{\m}(R)}$ is the unique maximal proper submodule $N$
with respect to the property
$$\a^{\rup{tq}} F^e(N) \subseteq N \text{ for  all } q = p^e,$$
where $F^e \colon H_{\m}^d(R) \to H_{\m}^d(R)$ is the $e$-times iterated
Frobenius induced on $H_{\m}^d(R)$.
\endproclaim

Now we generalize Theorem 3.4, which is the main theorem of Section 3, to
the case of rational coefficients.

\proclaim{Theorem 6.8}
Let $t \ge 0$ be a fixed rational number, $R$ a normal $\Q$-Gorenstein
local ring essentially of finite type over a field and let $\a$ be a
nonzero ideal. Assume that $\a \subseteq R$ is reduced from characteristic
zero to characteristic $p \gg 0$, together with a log resolution $f \colon
X \to Y= \Spec R$  of the ideal $\a$ such that $\a\OO_X = \OO_X(-Z)$. Then
$$\tau(\a^t) = H^0(X,\OO_X(\rup{K_{X/Y}-tZ})).$$
\endproclaim

\demo{Sketch of the proof}
This is also proved in a similar way as Theorem 3.4, so we just indicate
the points where some modification is needed in the following. First, we
note that Lemma 3.6 holds for rational coefficients without changing the
proof, i.e., $\tau((\a S)^t) \cap R = \tau(\a^t)$ under the assumption of
Lemma 3.6, and that Lemma 3.7
is already proved for rational coefficients. Hence we can use a canonical
covering of $R$ to reduce the proof of Theorem 6.8 to the quasi-Gorenstein
case. Then it is sufficient to prove the following generalization of
Theorem 3.9.

\proclaim{Theorem 6.9}
Let $t \ge 0$ be a fixed rational number, $(R,\m)$ a $d$-dimensional
normal local ring essentially of finite type over a field and let
$\a$ be a nonzero ideal. Assume that $\a \subseteq R$ is reduced from
characteristic zero to characteristic $p \gg 0$, together with a log
resolution $f \colon X \to Y = \Spec R$ of $\a$ such that $\a\OO_X =
\OO_X(-Z)$. Then
$$0^{*\a^t}_{H_{\m}^d(R)} =
\Ker\left( H^d_{\m}(R) \overset{\delta}\to\lra H^d_E(\OO_X(tZ))\right),$$
where $E$ is the closed fiber of $f$ and $\delta$ is an edge map as
in $3.8$.
\endproclaim

Here we note that the canonical dual of the sheaf $\OO_X(tZ) = \OO_X
(\rd{tZ})$ is $\omega_X(-\rd{tZ}) = \omega_X(\rup{-tZ})$, which is
isomorphic to $\OO_X(\rup{K_{X/R} - tZ})$ via $\omega_R \cong R$ if
$R$ is quasi-Gorenstein.

The inclusion $0^{*\a^t}_{H_{\m}^d(R)} \supseteq \Ker\delta$ of the
above theorem holds true in arbitrary fixed characteristic $p > 0$:
Just take an element $c$ in the proof of Proposition 3.8 from
$\a^{\rup{tq}} \subseteq H^0(X,\OO_X(-tqZ))$ instead of $\a^q =
H^0(X,\OO_X(-qZ))$, which gives rise to a map $cF^e \colon H^d_E
(\OO_X(tZ)) \to H^d_E(\OO_X(tZ))$. Then one sees that $\a^{\rup{tq}}
F^e(\Ker\delta) \subseteq \Ker\delta$ for all $q = p^e$, and
Proposition 6.7 applies.

To prove the reverse inclusion $0^{*\a^t}_{H_{\m}^d(R)} \subseteq
\Ker\delta$, we choose, in characteristic zero before reducing to
characteristic $p \gg 0$, a nonzero element $c \in \a$ such that
$R_c$ is regular and a log resolution $f \colon X \to \Spec R$ of
the ideal $c\a$, as in the proof of Theorem 3.9. Then choose an
$f$-ample $f$-exceptional $\Q$-Cartier divisor $D$ and a sufficiently
small $\varepsilon > 0$ so that $\rd{\tilde{Z} + \varepsilon\Div_X(c)}
= \rd{tZ}$, where $\tilde{Z}= tZ-D$. We then move to reduction modulo
$p \gg 0$ and let $m = r+2s+\rup{n_0t}$, keeping the integers $r,s,n_0$
just as same as in the proof of Theorem 3.9, i.e., the ideal $\a$ is
generated by $r$ elements, $c^s$ is a usual test element and also
an $\a^t$-test element, and $\Cal{K} = \bigoplus_{n>0} H^0(X,\omega_X
(-\rd{n\tilde{Z}}))$ is generated in degree $\le n_0$ as a graded
module over $\Cal{R} = \bigoplus_{n\ge 0} H^0(X,\OO_X(-n\tilde{Z}))$.
Now, arguing as before,
we obtain the required inclusion $0^{*\a^t}_{H_{\m}^d(R)} \subseteq
\Ker\delta$.
\enddemoo

\medskip
Finally, we note that the results in Section 4 also generalize to
rational coefficients, with the same proof; see \cite{DEL}, \cite{How},
\cite{La} for the corresponding results for multiplier ideals.

\proclaim{Theorem 6.10}
Let $t,t'$ be any nonnegative rational numbers.
\roster
\item {\rm (Restriction theorem):}
Under the assumption of Theorem $4.1$ we have
$$\tau((\a S)^t) \subseteq \tau(\a^t)S.$$

\item {\rm (Subadditivity in regular local rings;} cf.\ Remark $6.6)
\colon$ Under the assumption of Theorem $4.5$ we have
$$\tau(\a^t\bb^{t'}) \subseteq \tau(\a^t)\tau(\bb^{t'}).$$

\item Under the assumption of Theorem $4.8$, let $\a \subseteq R$ be
a toric ideal. Then $\tau(\a^t) = \J(\a^t)$, and it is also a toric ideal
generated by monomials $x^m$ with $m \in M$ such that
$$m + w \in \Int(t\cdot P(\a)).$$
\endroster
\endproclaim


\bigskip
Mathematical Institute, Tohoku University, Sendai 980--8578, Japan

{\it E-mail address}: {\tt hara\@math.tohoku.ac.jp}

\medskip
Graduate School of Mathematics, Nagoya University, Chikusa-ku,
Nagoya 464--8602, Japan

{\it E-mail address}: {\tt yoshida\@math.nagoya-u.ac.jp}


\Refs
\widestnumber\key{HWY2}

\ref\key{AKM}
\by I. Aberbach, M. Katzman, and B. MacCrimmon
\paper Weak F-regularity deforms in $\Q$-Gorenstein rings
\jour J. Algebra
\vol 204
\yr 1998
\page 281--285
\endref

\ref\key{AM}
\by I. Aberbach and B. MacCrimmon
\paper Some results on test ideals
\jour Proc. Edinburch Math. Soc. (2)
\vol 42
\yr 1999
\page 541--549
\endref

\ref\key{B}
\by J.-F. Boutot
\paper Singularit\'es rationelles et quotients par les groupes
r\'eductifs
\jour Invent. Math.
\vol 88 \yr 1987 \page 65--68
\endref

\ref\key{BS}
\by J. Brian\c con and H. Skoda
\paper Sur la cl$\hat{o}$ture int\'egrale d\'un id\'eal
de germes de fonctions holomorphes en un point de $C^n$
\jour C. R. Acad. Sci. Paris S\'er. A
\vol 278
\yr 1974
\page 949--951
\endref

\ref\key{BH}
\by W. Bruns and J. Herzog
\book  Cohen--Macaulay Rings
\bookinfo Cambridge Studies in Advanced Mathematics
\vol 39
\publ Cambridge University Press, Cambridge
\yr 1993
\endref

\ref\key{DEL}
\by J.-P. Demailly, L. Ein and R. Lazarsfeld
\paper A subadditivity property of multiplier ideals
\jour Michigan Math. J.
\vol 48
\yr 2000
\page 137--156
\endref

\ref\key{Ei}
\by L. Ein
\book Multiplier ideals, vanishing theorems and applications, {\rm in}
Algebraic Geometry---Santa Cruz 1995
\bookinfo pp. 203--219, Proc. Symp. Pure Math.
\vol 62
\yr 1997
\publ American Mathematical Society, Providence
\endref

\ref\key{ELS}
\by L. Ein, R. Lazarsfeld and K. E. Smith
\paper Uniform bounds and symbolic powers on smooth varieties
\jour Invent. Math.
\vol 144
\yr 2001
\page 241--252
\endref

\ref\key{FW}
\by R. Fedder and K.-i. Watanabe
\book  A characterization of F-regularity in terms of F-purity, {\rm in}
Commutative Algebra, Berkeley 1987
\bookinfo pp. 227--245, Math. Sci. Res. Inst. Publ.
\vol 15
\publ Springer-Verlag, New York
\yr 1989
\endref

\ref\key{GS}
\by S. Goto and Y. Shimoda
\book On the Rees algebras of Cohen-Macaulay local rings, {\rm in}
Commutative algebra, Fairfax 1979
\bookinfo pp. 201--231, Lecture Notes in Pure and Appl. Math.
\vol 68
\yr 1982
\publ Dekker, New York
\endref

\ref\key{GW}
\by S. Goto and K.-i. Watanabe
\paper On graded rings, I
\jour J. Math. Soc. Japan
\vol 30
\yr 1978
\pages 179--213
\endref

\ref\key{Ha1}
\by N. Hara
\paper A characterization of rational singularities in terms of
injectivity of Frobenius maps
\jour Amer. J. Math.
\vol 120
\yr 1998
\pages 981--996
\endref

\ref\key{Ha2}
\bysame
\paper Geometric interpretation of tight closure and test ideals
\jour Trans. Amer. Math. Soc.
\vol 353
\yr 2001
\page 1885--1906
\endref

\ref\key{HT}
\by N. Hara and S. Takagi
\paper Some remarks on a generalization of test ideals
\jour preprint
\endref

\ref\key{HW}
\by N. Hara and K.-i. Watanabe
\paper F-regular and F-pure rings vs.\ log terminal and
log canonical singularities
\jour J. Algebraic Geometry
\vol 11
\yr 2002
\page 363--392
\endref

\ref\key{HWY1}
\by N. Hara, K.-i. Watanabe, and K. Yoshida
\paper F-rationality of Rees algebras
\jour J. Algebra
\yr 2002
\vol 247
\page 153--190
\endref

\ref\key{HWY2}
\bysame
\paper Rees algebras of F-regular type
\jour J. Algebra
\yr 2002
\vol 247
\page 191--218
\endref

\ref\key{HHK}
\by M. Herrmann, E. Hyry and T. Korb
\paper On Rees algebras with a Gorenstein Veronese subring
\jour J. Algebra
\vol 200 \yr 1998 \pages 279--311
\endref

\ref\key{Ho1}
\by M. Hochster
\paper Cyclic purity versus purity in excellent Noetherian rings
\jour Trans. Amer. Math. Soc.
\vol 231
\yr 1977
\page 463--488
\endref

\ref\key{Ho2}
\by M. Hochster
\paper The tight integral closure of a set of ideals
\jour J. Algebra
\vol 230
\yr 2000
\pages 184--203
\endref

\ref\key{HH0}
\by M. Hochster and C. Huneke
\paper Tight Closure and strong F-regularity
\jour Mem. Soc. Math. France
\vol 38
\yr 1989
\pages 119--133
\endref

\ref\key{HH1}
\by M. Hochster and C. Huneke
\paper Tight Closure, invariant theory, and the Brian\c con-Skoda theorem
\jour J. Amer. Math. Soc.
\vol 3
\yr 1990
\pages 31--116
\endref

\ref\key{HH2}
\bysame
\paper F-regularity, test elements, and smooth base change
\jour Trans. Amer. Math. Soc.
\vol 346
\yr 1994
\pages 1--62
\endref

\ref\key{HH3}
\bysame
\paper Tight closure in equal characteristic zero
\jour to appear
\endref

\ref\key{How}
\by J. A. Howald
\paper Multiplier ideals of monomial ideals
\jour Trans. Amer. Math. Soc.
\vol 353
\yr 2001
\page 2665--2671
\endref

\ref\key{Hu}
\by C. Huneke
\book Tight Closure and Its Applications
\bookinfo C.B.M.S. Regional Conf. Ser. in Math. No. 88
\yr 1996
\publ American Mathematical Society, Providence
\endref

\ref\key{Hy1}
\by E. Hyry
\paper Blow-up rings and rational singularities
\jour Manuscripta Math.
\vol 98
\yr 1999
\pages 377--390
\endref

\ref\key{Hy2}
\bysame
\paper Coefficient ideals and the Cohen-Macaulay property
of Rees algebras
\jour Proc. Amer. Math. Soc.
\vol 129
\yr 2001
\page 1299--1308
\endref

\ref\key{Ka}
\by Y. Kawamata
\paper The cone of curves of algebraic varieties
\jour Ann. Math.
\vol 119
\yr 1984
\page 603--633
\endref

\ref\key{Ku}
\by E. Kunz
\paper On Noetherian rings of characteristic $p >0$
\jour Amer. J. Math.
\vol 98
\yr 1976
\pages 999--1013
\endref

\ref\key{La}
\by R. Lazarsfeld
\paper Positivity in Algebraic Geometry 
\jour preprint 
\endref

\ref\key{Li}
\by Lipman, J.
\paper Adjoints of ideals in regular local rings
\jour Math. Res. Letters
\vol 1
\yr 1994
\pages 739--755
\endref

\ref\key{Mc}
\by B. MacCrimmon
\paper Weak F-regularity is strong F-regularity for rings with isolated
non-$\Q$-Gorenstein points
\jour Trans. Amer. Math. Soc.
\page to appear
\endref

\ref\key{Ma}
\by H. Matsumura
\book Commutative ring theory
\bookinfo Cambridge Studies in Advanced Mathematics
\vol 8
\publ Cambridge University Press, Cambridge
\yr 1986
\endref

\ref\key{MS}
\by V. B. Mehta and V. Srinivas
\paper A characterization of rational singularities
\jour Asian J. Math.
\vol 1
\yr 1997
\page 249--278
\endref

\ref\key{N}
\by A. Nadel
\paper Multiplier ideal sheaves and K\"ahler-Einstein metrics of
positive scalar curvature
\jour Ann. Math.
\vol 132
\yr 549--596
\endref

\ref\key{Si}
\by A. K. Singh
\paper F-regularity does not deform
\jour Amer. J. Math.
\vol 121
\yr 1999
\page 919--929
\endref

\ref\key{Sm1}
\by K. E. Smith
\paper F-rational rings have rational singularities
\jour Amer. J. Math.
\vol 119
\yr 1997
\pages 159--180
\endref

\ref\key{Sm2}
\bysame
\paper The multiplier ideal is a universal test ideal
\jour Comm. Algebra 
\vol 28
\yr 2000
\page 5915--5929
\endref

\ref\key{T}
\by S. Takagi
\paper An interpretation of multiplier ideals via tight closure
\jour preprint
\endref

\ref\key{VV}
\by P. Valabrega and G. Valla
\paper Form rings and regular sequences
\jour Nagoya Math. J.
\vol 72
\yr 1978
\pages 93--101
\endref

\ref\key{Vr}
\by A. Vraciu
\paper Strong test ideals
\jour J. Pure Appl. Algebra
\vol 167
\yr 2002
\page 361--373
\endref

\ref\key{Wa}
\by K.-i. Watanabe
\paper F-regular and F-pure normal graded rings
\jour J. Pure Appl. Algebra
\vol 71
\page 341-350
\yr 1991
\endref

\ref\key{Wi}
\by L. J. Williams
\paper Uniform stability of kernels of Koszul cohomology indexed by the
Frobenius endomorphism
\jour J. Algebra
\vol 172
\page 721--743
\yr 1995
\endref

\endRefs
\enddocument